\documentclass{gthack}   
\usepackage[margin=1.4in, left=1.45in, right=1.2in]{geometry}
\usepackage{graphicx}  
\usepackage{enumitem}
\title[Quasipositive links and Stein surfaces]{Quasipositive links and Stein surfaces}

\author{Kyle Hayden}
\givenname{Kyle}
\surname{Hayden}
\address{Columbia University, New York, NY 10027}
\email{hayden@math.columbia.edu}


\arxivreference{1703.10150}
\arxivpassword{56w49}

\theoremstyle{plain}
\newtheorem{thm}{Theorem}[section]    
\newtheorem{lem}[thm]{Lemma}          
\newtheorem{prop}[thm]{Proposition}
\newtheorem{cor}[thm]{Corollary}

\theoremstyle{definition}
\newtheorem{defn}[thm]{Definition}    
\newtheorem{ex}[thm]{Example}
\newtheorem{rem}[thm]{Remark}             
\newcommand{\cc}{\mathbb{C}}
\newcommand{\rr}{\mathbb{R}}
\newcommand{\st}{{\mathrm{st}}}
\newcommand{\id}{\operatorname{id}}

\newcommand{\Symp}{\operatorname{Symp}}

\newcommand{\Diff}{\operatorname{Diff}}

\begin{document}

\begin{abstract}    
We study the generalization of quasipositive links from the three-sphere to arbitrary closed, orientable three-manifolds. Our main result shows that the boundary of any smooth, properly embedded complex curve in a Stein domain is a quasipositive link. This generalizes a result due to Boileau and Orevkov, and it provides the first half of a topological characterization of links in three-manifolds which bound complex curves in a Stein filling. Our arguments replace pseudoholomorphic curve techniques with a study of characteristic and open book foliations on surfaces in three- and four-manifolds.  
\end{abstract}

\maketitle


\section{Introduction} \label{sec:intro}

The tradition of studying complex varieties in $\cc^n$ via their intersection with smooth hypersurfaces has played a motivating role in knot theory almost since its inception. For $n=2$, the study of braids has been tied especially tightly to this line of inquiry. In \cite{rudolph:qp-alg}, Rudolph introduced the notion of quasipositive braids and links in $S^3$ and showed that each such link arises as the transverse intersection of the unit sphere $S^3 \subset \cc^2$ with a smooth algebraic curve $f^{-1}(0)\subset \cc^2$, where $f$ is a non-constant polynomial. More generally, Rudolph conjectured that any link in $S^3$ arising as the transverse intersection between a smooth complex curve and the unit sphere in $\cc^2$ is quasipositive. This was confirmed by Boileau and Orevkov in \cite{bo:qp}, whose proof uses Gromov's celebrated theory of pseudoholomorphic curves.

Many of these ingredients  have natural analogs when $\cc^2$ is replaced with an arbitrary Stein surface, a complex surface $X$ that admits a proper biholomorphic embedding as a closed subset in some $\cc^n$. In this setting, the role of the unit three-sphere in $\cc^2$ is played by a hypersurface $Y \subset X$ of constant radius in $\cc^n$. The region in $X$ bounded by $Y$ is a \emph{Stein domain}. Any such level set $Y$ in $X$ has a natural contact structure $\xi$ given by complex tangent lines, and we can define \emph{quasipositive} braids and links with respect to open book decompositions of $Y$ compatible with $\xi$; see Section~\ref{sec:qp}. Our main theorem says that the analog of Rudolph's conjecture also holds in this setting.

\begin{thm}\label{thm:qp}
The boundary of any smooth, properly embedded complex curve in a Stein domain is a quasipositive link.
\end{thm}

In particular, we show that any such link is transversely isotopic to the closure of a quasipositive braid with respect to an open book whose monodromy is a product of positive Dehn twists along homologically nontrivial simple closed curves; motivated by analogous terminology for Lefschetz fibrations, we call this a \emph{positive allowable} open book. 
In contrast with \cite{bo:qp}, our argument entirely avoids Gromov's theory of pseudoholomorphic curves and any results about positivity of intersections between such curves \cite{micallef-white,mcduff:positivity}. Instead, we study complex curves in Stein surfaces via their \emph{characteristic foliations} (described precisely in Definition~\ref{def:stein-char}). We reduce this to a study of (possibly singular) surfaces in three-manifolds, where we prove new results relating characteristic and open book foliations; see Theorems~\ref{thm:embedded} and \ref{thm:foln} below.

We also prove a partial converse to Theorem~\ref{thm:qp}.

\begin{thm}\label{thm:symp}
If $K \subset Y$ is the closure of a quasipositive braid with respect to a positive allowable open book, then $K$ bounds a properly embedded symplectic surface in a Stein filling of $Y$.\end{thm}

This offers evidence for a stronger converse to Theorem~\ref{thm:qp}: a symplectic surface can always be made pseudoholomorphic by a deformation of the ambient almost-complex structure compatible with the natural symplectic structure on $X$, though it is unclear whether or not this can be accomplished through a deformation of the Stein structure itself. In the case where $\partial X$ admits a planar open book decomposition, Hedden has shown in \cite{hedden:subcritical} that  quasipositivity is also a sufficient condition for a braid to bound a complex curve.  In a different direction, we point out that the relative version of the symplectic Thom conjecture implies that the surfaces constructed in the proof of Theorem~\ref{thm:symp} are genus-minimizing in their relative homology classes; see \cite{gk:caps}.

A few words on our notions of quasipositivity are in order. By a classical theorem of Alexander, every closed, orientable three-manifold admits an open book decomposition. Any link that is braided with respect to this open book can be represented in terms of a \emph{pointed monodromy} of a page of the open book decorated with marked points. We define quasipositive braids for any open book in terms of certain generators of the pointed mapping class group of the page. Moreover, given a closed, orientable contact three-manifold $(Y,\xi)$, we say that a transverse (resp.~smooth) link type $K$ is \emph{quasipositive} if there is an open book decomposition of $Y$ supporting $\xi$ in which $K$ is transversely (resp.~smoothly) isotopic to a quasipositive braid.  As evidenced by Theorems~\ref{thm:qp} and \ref{thm:symp}, such links exhibit close connections to complex and symplectic structures. 

We also study the natural generalization of strongly quasipositive links, a special class of nullhomologous quasipositive links. These links bound symplectic surfaces in a compact piece of the symplectization of $(Y,\xi)$, and these surfaces are in fact isotopic to Seifert surfaces in $Y$. When $(Y,\xi)$ is tight, these Seifert surfaces can be used to show that the Eliashberg-Bennequin inequality is sharp for all strongly quasipositive links. These same ideas can be used to generalize a result of Rudolph's \cite{rudolph:qp-obstruction}: the only strongly quasipositive knot in a symplectically fillable contact three-manifold that bounds a smooth, properly embedded ``slice'' disk in a symplectic filling is the unknot.

As mentioned above, a key step of our proof of Theorem~\ref{thm:qp} involves relating quasipositivity to properties of characteristic foliations on surfaces. These results are framed in terms of Morse functions on surfaces. Following \cite[\S9.3]{ec:book}, we say that a vector field $v$ on a surface $S$ is \emph{weakly gradient-like} for a function $f:S \to \rr$ if the zeroes of $v$ coincide with the critical points of $f$ and $df(v)>0$ outside the zeroes of $v$. By extension, we say that a characteristic or open book foliation on $S$ is \emph{weakly gradient-like} for $f: S \to \rr$ if it is directed by a vector field that is weakly gradient-like for $f$.

\begin{thm}\label{thm:embedded}
Let $(Y,\xi)$ be a contact three-manifold containing transverse links $K_0$ and $K_1$ that cobound an embedded surface $S \subset Y$, oriented so that $\partial S=K_1 \cup - K_0$.  If $K_0$ is transversely isotopic to a quasipositive braid with respect to a compatible open book $(B,\pi)$ and the characteristic foliation on $S$ contains only positive singular points and is weakly gradient-like for a Morse function on $S$, then $K_1$ is also transversely isotopic to a quasipositive braid with respect to the same open book.
\end{thm}

To prove Theorem~\ref{thm:embedded}, we strengthen the known relationship between characteristic foliations and open book foliations on surfaces in three-manifolds. Together with a result of Ito and Kawamuro \cite[Theorem 2.21]{ito-kawamuro:open}, the below theorem constitutes a partial analog of the Giroux correspondence in the setting of singular foliations on surfaces in contact three-manifolds.

\begin{thm}\label{thm:foln}
Let $(Y,\xi)$ be a contact three-manifold with compatible open book $(B,\pi)$, and suppose that $S \subset Y$ is an embedded surface whose characteristic foliation is weakly gradient-like with respect to a Morse function $f:S \to \rr$. After an isotopy of $(B,\pi)$ through open books supporting $\xi$, the open book foliation on $S$ may be assumed to be weakly gradient-like with respect to $f$ near all of its critical level sets.

Moreover, the isotopy of $(B,\pi)$ is induced by an  ambient isotopy that preserves $\xi$ outside an arbitrarily small neighborhood of each singular point of the characteristic foliation on $S$.
\end{thm}

In fact, we show that Theorems~\ref{thm:embedded} and \ref{thm:foln} hold not only for embedded surfaces $S \subset Y$ but for certain well-behaved singular surfaces. This will play a key role in our study of characteristic foliations on surfaces in four-manifolds in Section~\ref{sec:stein}.

\begin{rem}
After completing an early draft of this work, the author learned (\cite{hedden:email} and \cite{hedden-etnyre:email}) that some of the above results have also been obtained in work-in-progress by a collaboration between Inanc Baykur, John Etnyre, Matt Hedden, Keiko Kawamuro, and Jeremy Van Horn-Morris. By different means, they prove Theorem~\ref{thm:qp} for subcritical Stein fillings of $\#^k(S^1 \times S^2)$, viewed with a planar open book, as well as the same general results concerning quasipositive representatives of strongly quasipositive links in arbitrary open books and Theorem~\ref{thm:symp}. See \cite{ik:bennequin} for other recent work on strongly quasipositive braids in arbitrary open books.
\end{rem}

\noindent \textbf{Organization }  We begin by defining quasipositive links in three-manifolds in  Section~\ref{sec:qp}. The precise definition appears in Section~\ref{subsec:qp}, preceded by an extensive discussion of the framework of pointed open books we use to describe braids in arbitrary open books in Sections~\ref{subsec:contact-open}-\ref{subsec:pointed}. In Section~\ref{sec:foln}, we discuss characteristic and open book foliations and prove Theorems~\ref{thm:embedded} and \ref{thm:foln}.  Our study of Stein surfaces begins in Section~\ref{sec:stein}, where we obtain Theorem~\ref{thm:qp} as a consequence of stronger results for a class of surfaces that generalize complex curves; see Theorems~\ref{thm:full}-\ref{thm:prod}.  We close with the proof of Theorem~\ref{thm:symp} in Section~\ref{sec:con}.

\noindent \textbf{Acknowledgements } I am deeply grateful to John Baldwin, Eli Grigsby, and Matt Hedden for helpful discussions during the development of the ideas presented here. In addition, I am grateful to John Etnyre, Keiko Kawamuro, and Sebastian Baader for communications that improved this paper. Finally, I thank the anonymous referee for their many helpful comments. 

This work was partially supported by NSF grant DMS-1803584.

\section{Quasipositive links in three-manifolds}
\label{sec:qp}

In this section, we generalize the notion of quasipositivity to braids in arbitrary open books; see \cite{square,ik:bennequin} for related definitions. Our definitions are framed in the language of pointed open books  introduced in \cite{bvvv:equivalence}; we review and expand on this framework in Section~\ref{subsec:pointed}. Throughout this paper, we assume a familiarity with open book decompositions and contact geometry at the level of \cite{etnyre:obd}. For later reference, we collect a few of these ideas in Section~\ref{subsec:contact-open}.

\subsection{Open book decompositions and contact structures} \label{subsec:contact-open}

Recall that an open book decomposition of a three-manifold $Y$ is a pair $(B,\pi)$, where $B$ is an oriented, fibered link in $Y$ and $\pi: Y \setminus B \to S^1$ is a fibration. Moreover, we require that the closure $F_\theta$ of each fiber $\pi^{-1}(\theta)$ in $Y$ is a Seifert surface for $B$. The link $B$ is called the \emph{binding} of the open book, and the Seifert surfaces $F_\theta$ are the \emph{pages} of the open book. 

 We can also describe an open book decomposition abstractly using a pair $(F,\varphi)$, where $F$ is a compact surface homeomorphic to a page and $\varphi$ is a diffeomorphism of $F$ that fixes its boundary and encodes the monodromy of the fibration. In particular, this defines a closed three-manifold $Y_{(F,\varphi)}=(F\times I)/{\sim}_\varphi$ where the relation is defined by
\begin{align*}
(x,1)\sim_\varphi (\varphi(x),0),& \qquad x \in F\\
(x,t) \sim_\varphi (x,s),& \qquad x \in \partial F, \ \ t,s \in I.
\end{align*}
We say that $(F,\varphi)$ \emph{corresponds} to an open book decomposition $(B,\pi)$ of a three-manifold $Y$ if there is a diffeomorphism $Y_{(F,\varphi)}\to Y$ that pulls $\pi$ back to the natural fibration from ${Y_{(F,\varphi)}}\setminus \partial F$ to $S^1$.

An open book decomposition of a contact three-manifold $(Y,\xi)$ is said to be \emph{compatible} with the contact structure if $\xi$ has a contact 1-form $\alpha$ such that $d\alpha$ is an area form on each page and $\alpha$ is positive along the binding. This definition is motivated by a construction of Thurston-Winkelnkemper \cite{tw:existence} that associates a contact structure to each open book decomposition of a three-manifold. Celebrated work of Giroux \cite{giroux:geometry} extends this to a one-to-one correspondence between open books (up to positive Hopf stabilization) and contact structures (up to isotopy) on a three-manifold.

\subsection{Transverse links, braids, and pointed open books.} \label{subsec:pointed} A link  in a contact three-manifold $(Y,\xi)$ is called \emph{transverse} if it is everywhere transverse to $\xi$ and \emph{Legendrian} if it is everywhere tangent to $\xi$. All contact structures in this paper will be assumed to be co-oriented, and all oriented transverse links will be assumed to intersect the contact planes positively unless otherwise specified.

Given an open book decomposition on $Y$, an oriented link is said to be \emph{braided} if it is positively transverse to the pages of the open book. In the case where the open book decomposition is compatible with the contact structure, any braid is naturally braid isotopic to a transverse link. The ``transverse Alexander theorem'' asserts the converse. The general form of the theorem was proven by Pavelescu, generalizing an earlier result due to Bennequin for links in the three-sphere; see  \cite{pavelescu:braids} and \cite{bennequin}, respectively.

\begin{thm}[Pavelescu]\label{thm:alexander}
Given any compatible open book decomposition of a contact three-manifold $(Y,\xi)$, every transverse link is transversely isotopic to a braid.
\end{thm}

\begin{rem}\label{rem:arcs}
As discussed in \cite{pavelescu:braids},  Theorem~\ref{thm:alexander} holds more generally for oriented arcs $\ell$ in $Y$ that are positively transverse to $\xi$ and positively braided near their endpoints. Moreover, applying the ``Bennequin trick'' as articulated in \cite[Lemmas 1-2]{b-m:simplicity}, we see that any such arc $\ell$ within a transverse link $K$ may be braided via a contact isotopy that fixes $K \setminus \ell$.
\end{rem}

The classical  operations of Markov stabilization and destabilization have analogs for braids in an arbitrary open book; these will be defined precisely below. The following theorem  generalizes work of Wrinkle \cite{wrinkle} and Orevkov-Shevchishin \cite{os:markov}.

\begin{thm}[Pavelescu]\label{thm:markov}
Given any compatible open book decomposition of a contact three-manifold $(Y,\xi)$, two transverse braids are transversely isotopic if and only if they have a common positive stabilization.
\end{thm}

Following \cite{bvvv:equivalence}, we can  encode braids abstractly using \emph{pointed open book decompositions}. In particular, let $(F,\varphi)$ be an abstract open book corresponding to $(B,\pi)$ and let $P=\{p_1,\ldots,p_n\}$ be a set of distinct marked points in the interior of $F$. The monodromy $\varphi$ is isotopic to a diffeomorphism $\widehat \varphi$ of $F$ that fixes the marked points setwise and $\partial F$ pointwise, i.e. a diffeomorphism of the pair $(F,P)$ fixing $\partial F$.  The link in $Y_{(F,\widehat \varphi)}$ given by $(P\times I)/{\sim_{\widehat\varphi}}$  corresponds to a link $K$ in $Y$ that is braided with respect to $(B,\pi)$. We say that $K$ is \emph{encoded} by the pointed open book decomposition $(F,P,\widehat\varphi)$, and we refer to $\widehat \varphi$ as the \emph{pointed monodromy}. Unless otherwise specified, isotopy of the pointed monodromy will be through diffeomorphisms that fix both $P$ and $\partial F$ pointwise. An extension of the argument from \cite[\S4.4.2]{geiges:intro} shows that any braid $K$ in $(B,\pi)$ is encoded by an abstract open book.

\begin{figure}\center
\def\svgwidth{\linewidth}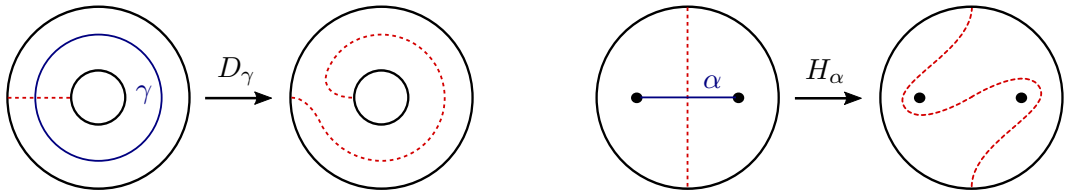 
\caption{A positive Dehn twist  \emph{(left)} and a positive half-twist \emph{(right)}.} \label{fig:arcNeighborhood}
\end{figure}

The pointed mapping class group of $F$ is generated by Dehn twists about simple closed curves and half-twists along embedded arcs joining distinct marked points; see Figure~\ref{fig:arcNeighborhood} for local depictions of these twists. We will let $D_\gamma$ denote a positive Dehn twist about a curve $\gamma$ and $H_\alpha$ denote a positive half-twist about an arc $\alpha$, respectively.  For later reference, we highlight the behavior of these twists  under conjugation:
\begin{lem}[{cf.~\cite[Fact 3.7]{farb-margalit}}]\label{lem:conjugation}
Let $\varphi$ be any diffeomorphism of a surface $F$.
\begin{enumerate}[label=(\alph*)]
\item For any simple closed curve $\gamma$ in $F$, we have $\varphi \circ D_\gamma \circ \varphi^{-1} = D_{\varphi(\gamma)}$.
\item For any embedded arc $\alpha$ in the interior of $F$, we have $\varphi \circ H_{\alpha} \circ \varphi^{-1} = H_{\varphi(\alpha)}$.
\end{enumerate}
\end{lem}

\begin{figure} \center
\def\svgwidth{340pt}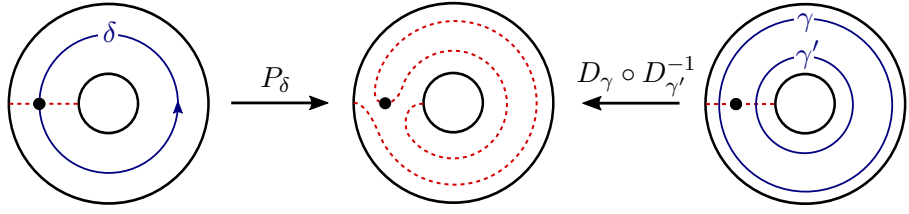 
\caption{Expressing a point-pushing map as a product of Dehn twists.}
 \label{fig:push}
\end{figure}

We also reserve the notation $P_\delta$ for the ``point-pushing'' map along a simple closed curve $\delta$ through a marked point; as shown in Figure~\ref{fig:push}, $P_\delta$ can be expressed as a product of oppositely-signed Dehn twists along curves parallel to $\delta$. The behavior of point-pushing maps under conjugation is analogous to that of Dehn twists and half-twists.

The next lemma characterizes braid isotopy in terms of abstract pointed open books.

\begin{lem}\label{lem:isotopy}
Suppose $(F,P,\widehat \varphi)$ and $(F,P',\widehat \varphi')$ are two pointed open books encoding a pair of braids $K$ and $K'$ in the same open book $(B,\pi)$, where $P$ and $P'$ are sets of $n$ marked points. The braids $K$ and $K'$ are braid isotopic if and only if $\widehat\varphi'$ is isotopic to $h \circ \widehat\varphi \circ h^{-1}$ for a diffeomorphism $h$ of $F$ that sends $P$ to $P'$ and is the identity near $\partial F$.
\end{lem}

\begin{proof}
We begin by noting that the braid isotopy extends to an ambient isotopy of $Y \setminus B$  that preserves each page of the open book, carries $K$ to $K'$, and is the identity near $B$. This gives rise to an isomorphism of the mapping tori of $\widehat \varphi \in \Diff(F,P)$ and $\widehat \varphi' \in \Diff(F,P')$. Any such isomorphism lifts to a diffeomorphism $f$ of $F \times [0,1]$ of the form $(x,t) \mapsto (f_t(x),t)$ for a smooth family $f_t \in \Diff(F)$. Moreover, because the braids in the mapping tori induced by $P$ and $P'$ (i.e.~$K$ and $K'$) are exchanged by the isomorphism, we have $f_t(P)=P'$ for all $t$.

In order to descend to the mapping torus of $\widehat \varphi$,  $f$ must satisfy $f(x,1)=f(\widehat \varphi(x),0)$, i.e.~$(f_1(x),1)=(f_0(\widehat \varphi(x)),0)$, for all $x \in F$. Similarly, since $(f_1(x),1)$ is identified with $(\widehat \varphi'(f_1(x)),0)$ in the mapping torus of $\widehat \varphi'$,  we must have $\widehat \varphi' \circ f_1 = f_0 \circ \widehat \varphi$, i.e.~$\widehat \varphi'=f_0 \circ \widehat \varphi \circ f_1^{-1}$. Now define $\widehat \varphi'_t= \widehat \varphi' \circ (f_1 \circ f^{-1}_{1-t})$. Then $\widehat \varphi'_t$ is an isotopy between $\widehat \varphi'_0=\widehat \varphi'$ and $\widehat \varphi'_1=f_0 \circ \widehat \varphi \circ f_0^{-1}$. Moreover, since $f_t(P)=P'$ for all $t$, we have $\widehat \varphi'_t \in \Diff(F,P')$ for all $t$. This shows that $\widehat \varphi'$ is isotopic in $\Diff(F,P')$ to a conjugate of $\widehat \varphi$, where the conjugating diffeomorphism carries $P$ to $P'$.
\end{proof}

We can also describe braid (de)stabilization in terms of abstract pointed open books: Let $p_{n+1}$ be a point in a collar neighborhood of $\partial F$ on which $\widehat\varphi$ is the identity and let $\alpha$ be an arc joining $p_{n+1}$ to $p_n$ (and meeting no other marked points). We say that the pointed open book $(F,P\cup\{p_{n+1}\}, H_\alpha \circ \widehat \varphi)$ is a \emph{positive Markov stabilization} of the original pointed open book. If $K$ and $K'$ in $Y$ are the corresponding braids with respect to the open book $(B,\pi)$, then $K'$ is obtained from $K$ by positive braid stabilization in the same sense defined by Pavelescu as referenced in Theorem~\ref{thm:markov}. Negative stabilizations and positive and negative destabilizations are defined similarly and also agree with the corresponding notions defined by Pavelescu.

The next lemma relates stabilization and the addition of half-twists.

\begin{lem}\label{lem:stab}
Suppose the braid $K_1$ is obtained from the braid $K_0$ by adding a half-twist. If $K_0^+$ is a positive Markov stabilization of  $K_0$, then $K_1$ has a positive Markov stabilization $K_1^+$ that is obtained from $K_0^+$ by adding a half-twist of the same sign.
\end{lem}

\begin{proof}Without loss of generality, we assume the half-twist is positive. Write the pointed monodromy for $K_0$ as $\widehat \varphi_0$ and the monodromy for $K_1$ as $\widehat \varphi_1= H_\alpha \circ \widehat \varphi_0$, where  $H_\alpha$ is the added half-twist.  If $K_0^+$ is obtained from $K_0$ by braid isotopy and a single positive stabilization, then it is encoded by a pointed open book $(F,P^+,\widehat \varphi_0^+)$ with  
$\widehat \varphi_0^+ = h \circ \widehat \varphi_0 \circ h^{-1} \circ H_{\beta},$
 where $h$ is a diffeomorphism of $F$ fixing a neighborhood of the boundary, $P^+$ is the union of $h(P)$ and a new marked point $p^+_{n+1}$ near the boundary, and $\beta$ is an arc joining $p^+_{n}=h(p_n)$ to $p^+_{n+1}$. Define a new braid $K_1^+$ via the pointed open book $(F,P^+,\widehat \varphi_1^+)$ where
 $$\widehat \varphi_1^+=h \circ H_\alpha \circ \widehat \varphi_0 \circ h^{-1} \circ H_\beta= H_{h(\alpha)} \circ h \circ \widehat \varphi_0 \circ h^{-1} \circ H_\beta = H_{h(\alpha)} \circ \widehat\varphi_0^+.$$
The first equality shows that $K_1^+$ is obtained from $K_1$ by braid isotopy and a single positive stabilization, and the last equality shows that $K_1^+$ is obtained from $K_0^+$ by adding a positive half-twist, as desired.

In general, if $K_0^+$ is obtained from $K_0$ by applying a sequence of braid isotopies and multiple Markov stabilizations, we may make an analogous argument by induction on the number of Markov stabilizations in the sequence.
\end{proof}

Finally, we consider the operation of Hopf stabilization on pointed open books. Suppose the transverse link $K$ in $(Y,\xi)$ is encoded by a pointed open book $(F,P,\widehat\varphi)$. An abstract positive Hopf stabilization of this pointed open book is a triple $(F', P, \widehat \varphi \circ D_\gamma)$, where $F'$ is obtained from $F$ by attaching an oriented 1-handle and $D_\gamma$ is a positive Dehn twist around a curve $\gamma \subset F'$ passing through this 1-handle exactly once and avoiding the marked points. This operation preserves the (transverse) link type:

\begin{lem}[{cf.~\cite[Corollary 2.5]{bvvv:equivalence}}]
With notation as above, the pointed open books $(F,P,\widehat\varphi)$ and $(F',P,\widehat \varphi \circ D_\gamma)$ represent equivalent (transverse) links in $(Y,\xi)$.
\end{lem}

\begin{proof}
Positive Hopf stabilization preserves the underlying contact three-manifold and simply changes the open book structure in a standard Darboux ball $N$ with convex boundary. In particular, $N$ is a small neighborhood of an arc in $Y$ that lies in a page $F_\theta$ of the original open book and is identified with the properly embedded arc $\gamma \cap F$; see \cite[\S2]{etnyre:obd} and \cite{gabai:murasugi}.  (This can also be seen explicitly by viewing the original and stabilized open books as two different ways of extending a partial open book decomposition for $(Y_0=Y \setminus \mathring{N},\xi|_{Y_0})$ to an open book for $(Y,\xi)$; see \cite{hkm:sutured} for background on partial open book decompositions.)  Moreover, since $\gamma$ avoids the set of marked points $P$, the braids represented by $(F,P,\widehat\varphi)$ and $(F',P,\widehat \varphi \circ D_\gamma)$ lie outside $N$; here the open book structures coincide. It follows that these braids are naturally identified and, in particular, represent the same transverse link in $(Y,\xi)$.
\end{proof}

\subsection{Quasipositive braids in an open book}\label{subsec:qp}

In this section we consider open book decompositions of arbitrary closed, oriented three-manifolds. We begin with the most basic notion of quasipositivity.

\begin{defn}\label{def:qp}
A braid in an  open book is \emph{quasipositive} if it can be encoded by an abstract pointed open book whose pointed monodromy is isotopic to a product of positive half-twists and arbitrary Dehn twists.
\end{defn}

\begin{rem}
To supply more context for this definition, we adopt the perspective of Etnyre-Van Horn-Morris in \cite{E-VHM:monoids}, which views classical quasipositive braids as the monoid generated by positive half-twists in the pointed mapping class group of the disk. We can say something similar in the more general case. Given a set of marked points $P$ in $F$, the elements of the pointed mapping class group that can be expressed as products of positive half-twists and arbitrary Dehn twists form a monoid. By Lemma~\ref{lem:conjugation}, it is easy to see that this monoid is preserved by conjugation. In particular, this shows that Definition~\ref{def:qp} is independent of the choice of the abstract pointed open book $(F,P,\widehat \varphi)$ used to encode a braid $K$ in $(B,\pi)$. Indeed, if $K$ is also encoded by $(F,P',\widehat\varphi')$ and $(F,P,\widehat\varphi)$ is quasipositive, then $\widehat \varphi'$ is isotopic to $h \circ \widehat \varphi \circ h^{-1}$ for some diffeomorphism $h$ of $F$ (as per the discussion in Section~\ref{subsec:pointed}) and is therefore quasipositive.
\end{rem}

We also generalize Rudolph's notion of strongly quasipositive braids \cite{rudolph:kauffman-bound}:

\begin{defn}\label{def:sqp}
A braid in an open book is \emph{strongly quasipositive} if it can be encoded by an abstract pointed open book $(F,P,\widehat\varphi)$ where the marked points $P$ lie along the interior edge of a collar neighborhood of $\partial F$ and $\widehat\varphi$ is a product of positive half-twists and arbitrary Dehn twists along arcs and curves whose interiors do not meet the collar neighborhood; see Figure~\ref{fig:sqp}.

\begin{figure} \center
\def\svgwidth{380pt}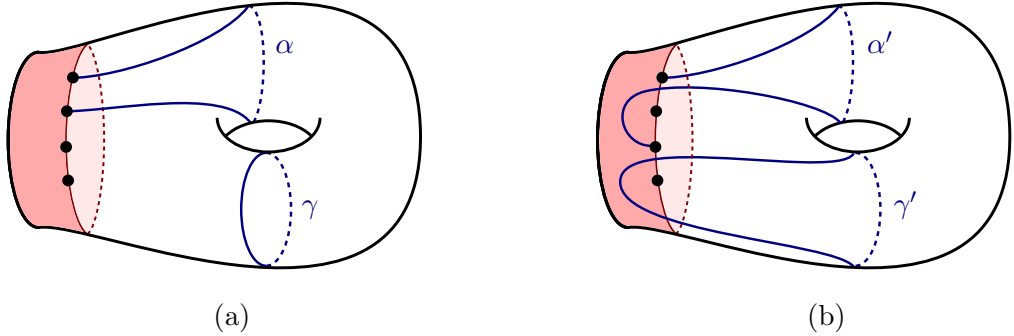 
\caption{Twists along $\alpha$ and $\gamma$ are allowed in a strongly quasipositive pointed monodromy, but twists along $\alpha'$ and $\gamma'$ are not.}
 \label{fig:sqp}
\end{figure}
\end{defn}

\begin{rem} 
Ito and Kawamuro present an equivalent definition  in \cite[Definition~4.3]{ik:bennequin}.
\end{rem}

For yet another perspective, we may associate to any marked surface $(F,P)$ a braid group $B_n(F,P)$, where $n=|P|$: the group $B_n(F,P)$ is defined to be the kernel of the forgetful map from the pointed mapping class group to the regular mapping class group of the underlying surface. Using Lemma~\ref{lem:conjugation} and the discussion that follows it, we see that the pointed monodromy of any pointed open book $(F,P,\widehat \varphi)$ can be separated into the composition of an element of the braid group $B_n(F,P)$ (expressed as a product of half-twists $H_{\alpha_1}^{\epsilon_1}\circ \cdots \circ H_{\alpha_\ell}^{\epsilon_\ell}$ with $\epsilon_i=\pm 1$) and a diffeomorphism of $F$ that fixes $P$ pointwise (expressed as a product of Dehn twists $D_{\gamma_1}^{\delta_1} \circ \cdots \circ D_{\gamma_k}^{\delta_k}$ with $\delta_i = \pm1$ along simple closed curves $\gamma_i$ missing $P$). We say that the braid encoded by $(F,P,\widehat \varphi)$ is the \emph{closure} of $H_{\alpha_1}^{\epsilon_1}\circ \cdots \circ H_{\alpha_\ell}^{\epsilon_\ell}  \in B_n(F,P)$ with respect to the  abstract open book $(F,D_{\gamma_1}^{\delta_1}\circ \cdots \circ D_{\gamma_k}^{\delta_k})$. This gives a natural refinement of the notion of quasipositivity with respect to a given factorization of the monodromy map for the underlying open book:

\begin{defn}\label{def:closure} An element of the braid group $B_n(F,P)$ is \emph{quasipositive} if it can be represented as a product of positive half-twists. A braid is \emph{quasipositive with respect to an (abstract) open book $(F,\varphi)$} if it can be represented as the closure of a quasipositive braid in $B_n(F,P)$ with respect to $(F,\varphi)$, where $\varphi$ fixes $P \subset F$ pointwise.
\end{defn}

As we will demonstrate in the proof of Theorem~\ref{thm:qp},  the boundary of a complex curve in a Stein domain is isotopic to the closure of a quasipositive braid with respect to a \emph{positive allowable} open book, i.e.~one where the monodromy for the underlying open book is a product of positive Dehn twists along homologically nontrivial curves.

Finally, as indicated in the introduction, we extend these definitions to transverse and smooth link types: Given a closed, orientable contact three-manifold $(Y,\xi)$, we say a transverse or smooth link type $K$ is \emph{(strongly) quasipositive} if there is an open book decomposition of $Y$ compatible with $\xi$ in which $K$ is transversely or smoothly isotopic, respectively,  to a (strongly) quasipositive braid.

\section{Characteristic and open book foliations}\label{sec:foln}

Contact structures and open books on a three-manifold $Y$ naturally induce singular foliations on surfaces $S \subset Y$. For general background on characteristic and open book foliations, we refer the reader to \cite{geiges:intro} and \cite{ito-kawamuro:open}, respectively. 

For a rigorous treatment of singular foliations, see \cite{geiges:intro} (especially Remark~2.5.21). Formally, we define an oriented, singular foliation $\mathcal{F}$ on a surface $S$ to be an equivalence class of vector fields $[v]$, where $[v]=[v']$ if there is a nowhere-zero function $f$ on $S$ such that $v'=fv$. However, we will often conflate a singular foliation with the associated decomposition of $S$ into integral curves of a representative vector field $v$ of $\mathcal{F}$. We will be concerned with singular foliations up to topological conjugacy; two singular foliations $\mathcal{F}$ and $\mathcal{F}'$ on a surface $S$ are said to be \emph{topologically conjugate} if there exists a homeomorphism of $S$ carrying $\mathcal{F}$ to $\mathcal{F}'$.

\begin{defn}\label{def:char}
Let $Y$ be a three-manifold containing a smooth, oriented surface $S$ and let $\Omega$ be an area form on $S$.
 The \emph{characteristic foliation} on $S \subset Y$ induced by $\xi=\ker \alpha$ is the singular foliation on $S$ directed by a vector field $v$ satisfying $\iota_v \Omega = \alpha|_S$.
\end{defn}
 
The existence of such a vector field $v$ is guaranteed by the nondegeneracy of the 2-form $\Omega$. It is straightforward to show that the oriented singular foliation on $S$ defined by $v$ depends only on $\xi=\ker \alpha$ and the co-orientation induced by $\alpha$ (cf \cite[\S2.5.4]{geiges:intro}).

An open book decomposition $(B,\pi)$ of $Y$ also defines a singular foliation on $S \subset Y$; see \cite{ito-kawamuro:open}.  This singular foliation admits a description analogous to Definition~\ref{def:char}. To set this up, let $d\theta$ denote the angular coordinate 1-form on $S^1$. Define a 1-form $\lambda=f  \pi^*(d\theta)$ on $Y$, where $f$ is strictly positive on $Y \setminus B$ and is locally given by $f(s,r,\vartheta)=r^2$ in some coordinates $(s,r,\vartheta) \in S^1 \times D^2$ near each component of $B$. When restricted to $Y \setminus B$,  $\ker \lambda$  is a co-oriented, integrable plane field whose maximal integral submanifolds are precisely the pages of the open book $(B,\pi)$.

\begin{defn}\label{def:open}
The \emph{open book foliation} on $S$ induced by $(B,\pi)$ is the singular foliation directed by a vector field $v$ satisfying $\iota_v \Omega =\lambda|_S$, with $\lambda$ as above.
\end{defn}

It is straightforward to verify that Definition~\ref{def:open}  coincides with the definition given in \cite{ito-kawamuro:open}. In particular, note that the regular leaves of the open book foliation on $S$ are tangent to the pages of $(B,\pi)$. 

These definitions of characteristic and open book foliations naturally extend to oriented surfaces $S$ equipped only with a smooth map $\phi: S \to Y$, where the conditions $\iota_v \Omega=\alpha|_S$ and $\iota_v \Omega= \lambda|_S$ from Definitions~\ref{def:char} and \ref{def:open} are replaced with the conditions $\iota_v \Omega= \phi^*\alpha$ and $\iota_v \Omega=\phi^* \lambda$. We will be particularly interested in maps $\phi: S \to Y$ that are well-behaved with respect to a Morse function on $S$.

\begin{defn}\label{def:nonsing}
Given a  map $\phi: S \to Y$ and a  function $f: S \to \rr$, we say that $\phi$ is \emph{nonsingular with respect to $f$} if $\phi$ restricts to an embedding on each level set of $f$ and on a neighborhood of each critical point of $f$.
\end{defn}

This definition is motivated by the following source of examples.

\begin{lem}\label{lem:nonsing}
Given an embedded surface $S \subset Y \times \rr$, the map $\phi: S \to Y$ (induced by $Y \times \rr \to Y$) is nonsingular with respect to $f: S\to \rr$ (induced by $Y \times \rr \to \rr$).
\end{lem}

\begin{proof}
The projection $Y\times \rr \to Y$ restricts to a diffeomorphism between each level set $Y_t=Y \times \{t\}$ and $Y$. Thus $\phi$ restricts to an embedding on each level set $f^{-1}(t)= S \cap Y_t$. Next, observe that the map $d\phi$ is nondegenerate wherever $S$ is transverse to the fibers $\{y\} \times \rr$ of the projection $Y\times \rr \to Y$. At any critical point $p$ of $f$, the surface $S$ is tangent to a level set $Y_t$, hence transverse to the fibers $\{y\} \times \rr$. It follows that $d\phi$ is nondegenerate at $p$, so $\phi$ is an immersion near $p$ and thus restricts to an embedding on a sufficiently small neighborhood of $p \in S$.
\end{proof}

Note that the level sets of $f$ are oriented by a vector field $w$ satisfying $\iota_w \Omega= -df$. (This coincides with the ``outward-normal first'' orientation convention for $f^{-1}(t) = \partial f^{-1}(-\infty,t]$.)

\subsection{Model subsurfaces}

Recall the definition of weakly gradient-like vector fields and singular foliations from Section~\ref{sec:intro}. Observe that if the characteristic or open book foliation on $S\subset Y$ is weakly gradient-like for $f: S \to \rr$, then the regular level sets of $f$ are transverse to the contact planes or pages of the open book, respectively. (This is true for any  smooth map $\phi: S \to Y$ that is nonsingular with respect to $f$.)

\begin{figure}\center
\def\svgwidth{\linewidth}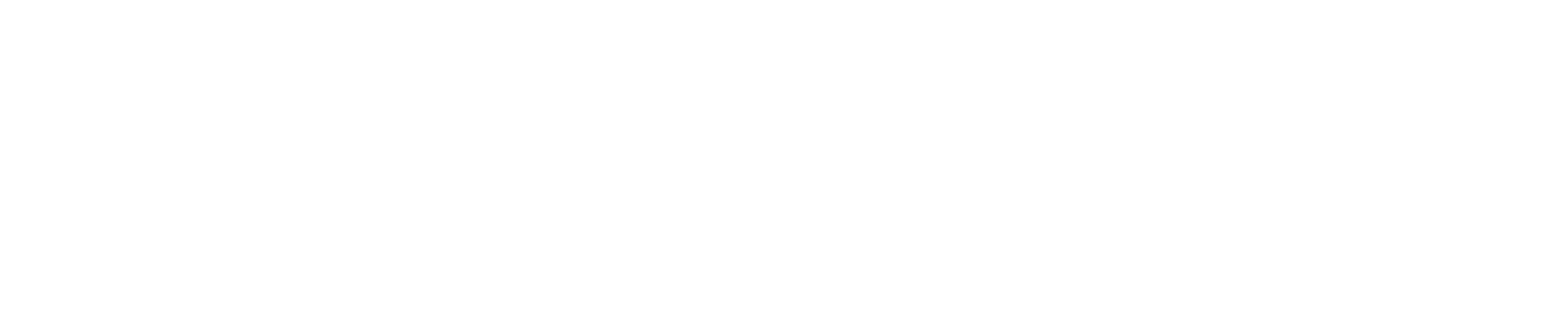 
\caption{A disk with elliptic foliation and a pair of pants with hyperbolic foliation are shown in (a) and (b), respectively. In (c), the boundary components of the pair of pants are contracted to cobound a band.}
 \label{fig:regions}
\end{figure}

 We will extract information from singular foliations on surfaces by decomposing them into  ``model subsurfaces'' that carry simple geometric and topological information. Recall that if $S$ is a disk or a pair of pants, then it admits an obvious Morse function $f$ with a single critical point of even or odd index, respectively. The singular foliation on $S$ induced by the gradient vector field of $f$ is unique up to topological conjugacy (and a choice of orientation).  See Figure~\ref{fig:regions}(a,b) for depictions of these models, which we refer to as a \emph{disk with elliptic foliation} and a \emph{pair of pants with hyperbolic foliation}, respectively. 

 We summarize the significance of these models in the context of  open book foliations:

\begin{lem}\label{lem:model}  Consider a three-manifold with a fixed open book decomposition.
\begin{enumerate}[label=(\alph*)]

\item The boundary of a disk with positive (resp. negative) elliptic open book foliation is a positively (resp. negatively) braided unknot meeting each page in one point.

\item If braids $K$ and $K'$ cobound a pair of pants with hyperbolic open book foliation (oriented so that its boundary is $K' \cup -K$), then $K'$ is obtained from $K$ by braid isotopy and the addition of a half-twist whose sign agrees with the sign of the hyperbolic point.
\end{enumerate}
Moreover, this holds for singular disks and pairs of pants, provided the maps are nonsingular with respect to the standard Morse functions described above.
\end{lem}

\begin{proof} The first  claim follows immediately from the standard definition of the open book foliation, as elliptic points correspond to intersections with the binding; see \cite{ito-kawamuro:open}. The second claim also follows quickly from definitions, following the method of proof in \cite[Lemma 1.6]{birman-finkelstein}. In particular, note that any isotopy of $K$ and $K'$ through curves transverse to the open book foliation on $S$ is a braid isotopy. Thus we may contract $K$ and $K'$ towards one another along the leaves of the open book foliation on $S$ until they coincide outside a neighborhood of the hyperbolic point, as depicted in Figure~\ref{fig:regions}(c). The region they enclose is a twisted band, and a local depiction of its embedding is given in Figure~\ref{fig:bennequin}(c) (which may be compared to \cite[Figures 1.12-1.13]{birman-finkelstein}). The claim follows.

In the case of a singular disk that is nonsingular with respect to  the standard Morse function, part (a) is immediate because any sufficiently small neighborhood of the elliptic point is embedded and  has an elliptic open book foliation. For part (b), let $f: P \to \rr$ denote the standard Morse function on the pair of pants $P$ with respect to which the open book foliation is weakly gradient-like and with respect to which the map $P \to Y$ is nonsingular. We may again contract each of $K$ and $K'$ along the leaves of the open book foliation until they both coincide with the critical level set of $f$ outside a neighborhood of the hyperbolic point. In particular, we may define this isotopy so that each intermediate link $K_t$ and $K_t'$ coincides with a level set of $f$ outside a fixed, embedded neighborhood of the hyperbolic point. This ensures that each intermediate link is embedded and braided. As above we may then locate an embedded, twisted band inside this neighborhood that these modified braids cobound.
\end{proof}

\subsection{Modifying open book foliations}
\label{subsec:modify}

We now prove Theorem~\ref{thm:foln}, which allows us to modify an open book on $(Y,\xi)$ so that the open book foliation on a surface $S$  approximates its characteristic foliation.

\begin{rem}\label{rem:embedded-nonsing} This theorem holds for a singular surface $S$ that is nonsingular with respect to the Morse function from the statement of the theorem, provided the collection of critical level sets is embedded. The proof given below carries over verbatim.
\end{rem}

\begin{proof}[Proof of Theorem~\ref{thm:foln}] We will prove something slightly stronger: Let $\Gamma$ be any finite collection of level sets of $f:S \to \rr$. We will produce an isotopy of any compatible open book $(B,\pi)$ through open books supporting $\xi$ so that the open book foliation on $S$ is weakly gradient-like with respect to $f$ on a neighborhood of $\Gamma$. For convenience, we assume the critical points of $f$ have distinct values; the general case follows with minor modifications.

We begin with the case of a regular level set $\gamma \subset  \Gamma$. The level set $\gamma$ is transverse to the characteristic foliation on $S$, so it is  a transverse link in $(Y,\xi)$. By the transverse Alexander theorem (Theorem~\ref{thm:alexander} and Remark~\ref{rem:arcs}), we may braid $\gamma$ using an ambient contact isotopy. Moreover, the isotopy may be taken to fix all other level sets in $\Gamma$. Equivalently, we may isotope $(B,\pi)$ by the inverse of this ambient isotopy while leaving $S$ fixed. After this isotopy of $(B,\pi)$, $\gamma$ is braided, and the same is true of all other level sets sufficiently close to $\gamma$.  It follows that the open book foliation is weakly gradient-like with respect to $f$ near $\gamma$.

Next consider a critical level set in $\Gamma$, which we have assumed contains a unique critical point. Near each regular component of the critical level set, we may apply the above argument. It thus suffices to consider a critical level set with no regular components. 

In the case of a minimum or maximum of $f$, i.e.~a positive or negative elliptic point of the characteristic foliation, the singular component of the level set is just the critical point $p \in S$ itself. Without loss of generality, we assume  $p$ is a  positive elliptic point (i.e.~local minimum) so that $T_p S=\xi_p$ with orientation. We can find (1) a path in $Y$ from  $p$ to the binding $B$ that meets $\Gamma$ only at $p$ and (2) an ambient contact isotopy that moves $p$ to the binding and is the identity outside a small neighborhood of the path \cite[Corollary~2.6.3]{geiges:intro}. This preserves the characteristic foliation on $S$ but turns $p$ into an elliptic point of the open book foliation (because $T_p S =\xi_p$ is transverse to the binding $B$).  Moreover, the open book foliation is transverse to all regular level sets of $f$ near $p$. Indeed, for some choice of local coordinates $(x,y)$ on $S$ centered at $p$, the open book foliation is directed by the radial vector field $v=x \partial_x+y \partial_y$.  It is easy to verify that $df(v)>0$ except at $p=(0,0)$, as desired, using the fact that $f$ has a local minimum at $p$. (For example, use these local coordinates to write $f(x,y)=f_{xx}(0,0)x^2 + f_{xy}(0,0)xy +f_{yy}(0,0) y^2+g(x,y)$, where  $g$ is a smooth function whose first and second derivatives all vanish at the origin. Then calculate the determinant of the Hessian for $df(v)$ using that of $f$.) As above, we may then pull $(B,\pi)$ back by the ambient contact isotopy to obtain the desired open book foliation on $S$  while leaving $S$ itself fixed.

Next suppose that the critical point $p$ is a saddle of $f$, ie~a hyperbolic singular point of the characteristic foliation. Without loss of generality, we assume $p$ is a positive hyperbolic point so that $T_p S=\xi_p$ with orientation. We will define a two-dimensional foliation $\mathcal{F}$ on a neighborhood of $p$ in $Y$ that locally induces a singular foliation on $S$ with a hyperbolic singularity; we will then perturb the open book $(B,\pi)$ so that its pages coincide with $\mathcal{F}$. To arrange this, fix a neighborhood $U$ of $p$ in $S$ and a neighborhood $U \times (-\epsilon,\epsilon)$ of $p$ in $Y$ in which $\xi$ is vertically invariant. By the Morse lemma, we may find coordinates $(x,y)$ on $U$ in which $f(x,y)=x^2-y^2$. Letting $z$ give the coordinate on $(-\epsilon,\epsilon)$, define a smooth 1-form by $\eta=dz+y \, dx + x \, dy$. Since $d\eta=0$, the plane field $\ker \eta$ is integrable and thus defines a two-dimensional foliation $\mathcal{F}$ on $U \times (-\epsilon,\epsilon)$ by Frobenius' theorem. Observe that the singular foliation on $S=U \times \{0\}$ is directed by the vector field $v=x \partial_x - y \partial_y$ (seen by noting that $v$ vanishes only at the origin and satisfies $\eta(v)=0$), which is gradient-like for $f$ near $p=(0,0,0)$. 

Let $R_\alpha$ be any Reeb vector field for $\xi$ that is positively transverse to the pages of $(B,\pi)$. We claim that $R_\alpha$ is also positive transverse to $\mathcal{F}$ near $p$. In the local coordinates above, $\xi$ coincides with $\ker dz$ at $p=(0,0,0)$. Thus we have $dz(R_\alpha)>0$ at $p=(0,0,0)$. Since $\eta_p=dz_p$, we have $\eta_p (R_\alpha) >0$ at $p$, hence $\eta(R_\alpha)>0$ in a neighborhood of $p$.

To set up the perturbation of $(B,\pi)$, let $\mathcal{F}_0$ denote the two-dimensional foliation of the neighborhood of $p$ induced by the pages of $(B,\pi)$. By restricting our neighborhood, we may assume that $R_\alpha$ is positively transverse to both $\mathcal{F}$ and $\mathcal{F}_0$. It follows that we can isotope $\mathcal{F}_0$ through a family of foliations $\mathcal{F}_t$ transverse to $R_\alpha$ such that $\mathcal{F}_t$ agrees with $\mathcal{F}_0$ outside some neighborhood of $p$ for all $t$ and $\mathcal{F}_1$ agrees with $\mathcal{F}$ inside some smaller neighborhood of $p$. (This can be achieved by a flow supported in the neighborhood that preserves the flowlines of $R_\alpha$.)  The pages remain positively transverse to $R_\alpha$, so this deformation defines a family of open books $(B,\pi_t)$ that support $\xi$. By construction,  the open book $(B,\pi_1)$ induces an open book foliation on $S$ with a hyperbolic point at $p$. Moreover, since $df(v) \geq 0$ with equality only at $p$, we see that the new open book foliation is weakly gradient-like with respect to $f$ in the neighborhood. 

We now extend this to the rest of the critical component $\gamma$. This level set $\gamma$  is a wedge of two circles; the subset $\gamma \setminus \{p\}$ is a pair of (oriented) open arcs that are transverse to $\xi$ and, near $p$, are braided with respect to $(B,\pi_1)$. Applying Theorem~\ref{thm:alexander} and Remark~\ref{rem:arcs} again, we may  braid $\gamma \setminus \{p\}$ with respect to $(B,\pi_1)$ via an ambient contact isotopy that fixes a neighborhood of $p$ and all other components of $\Gamma$. The open book foliation is preserved near $p$, and the portion of $\gamma$ lying outside the neighborhood of $p$ (as well as all nearby level sets) is braided with respect to the open book. It follows that the open book foliation is now gradient-like with respect to $f$ near $\gamma$. As above, we then pull the open book back by the isotopy so that we achieve the desired open book foliation on $S$ while leaving $S$ itself fixed.

This argument can be applied to each level set in $\Gamma \subset S$ without modifying the open book foliation near the rest of $\Gamma$, so the claim follows by induction on the size of $\Gamma$.
\end{proof}

\subsection{Characteristic foliations and quasipositivity}

We now prove Theorem~\ref{thm:embedded}. 

\begin{rem}\label{rem:foln-nonsing}
As with Theorem~\ref{thm:foln} (per Remark~\ref{rem:embedded-nonsing}), the proof of Theorem~\ref{thm:embedded} we give below holds for well-behaved singular surfaces in $Y$. In particular, $\phi: S \to Y$ must be nonsingular with respect to a Morse function $f: S \to \rr$ whose collection of critical level sets is embedded in $Y$, and the characteristic foliation on $S$ induced by $\phi$ must be weakly gradient-like for $f$ and contain only positive singular points.
\end{rem}

\begin{proof}[Proof of Theorem~\ref{thm:embedded}] We recall the key hypotheses: the braid $K_0$ is transversely isotopic to a  braid $K_0^+$ in an open book $(B,\pi)$ encoded by a quasipositive pointed open book; and the characteristic foliation on $S$ is weakly gradient-like with respect to a Morse function $f: S \to \rr$ with distinct critical values. We may assume that $K_0^+$ has been chosen to avoid a neighborhood of each critical point. (This can always be achieved by an arbitrarily small transverse braid isotopy of $K_0^+$, and it is easy to see that the hypotheses on the presentation of $K_0^+$ as a quasipositive braid are not affected by a small perturbation.) This simplifies a technical point: To apply Theorem~\ref{thm:foln}, we must modify $(B,\pi)$ using an ambient isotopy that restricts to a contact isotopy outside  an arbitrarily small neighborhood of each index-one critical point of $f: S \to \rr$. Since $K_0^+$ can be chosen to avoid these regions, it is carried along by a transverse isotopy to a braid in the modified open book $(B',\pi')$. Furthermore, we may use the ambient isotopy to carry forward the quasipositive pointed open book encoding $K_0^+$ in $(B,\pi)$ to a quasipositive pointed open book encoding $K_0^+$ in $(B',\pi')$.

Applying Theorem~\ref{thm:foln}, we may assume that the boundary $\partial S=K_1 \cup - K_0$ is braided with respect to $(B',\pi')$ and the open book foliation on $S$ is weakly gradient-like for $f$ near its critical level sets.

Divide $S$ into subsurfaces by splitting it along small neighborhoods $f^{-1}([c-\epsilon,c+\epsilon])$ of its critical level sets $f^{-1}(c)$. For sufficiently small $\epsilon>0$, the open book foliation on the surface $f^{-1}([c-\epsilon,c+\epsilon])$ is weakly gradient-like for $f$. Note this implies that each such subsurface has braided boundary. Treating each subsurface separately and arguing inductively on the number of critical points, it suffices to consider these cases: 
\begin{enumerate}[label=(\roman*)]
\item $S$ contains no critical points of $f$, or
\item the open book foliation on $S$ is weakly gradient-like for $f$ and $S$ has a unique critical of (a) even index, or (b) odd index.
\end{enumerate}
In case (i), the surface is a collection of annuli. Since $f$ has no critical points, the characteristic foliation on each annulus is linear. This implies that $K_1$ is transversely isotopic to $K_0$, which in turn is transversely isotopic to a quasipositive braid with respect to $(B',\pi')$.

In case (ii-a), the (unique) critical point of even index must correspond to a positive singular point of the characteristic foliation. The surface then consists of a collection of annuli with linear open book foliation and a disk with positive elliptic open book foliation, as depicted in Figure~\ref{fig:regions}. By Lemma~\ref{lem:model}, $K_1$ is obtained from $K_0$ by braid isotopy and the addition of a disjoint braided unknot $U$ that bounds a meridional disk of the binding. Since $K_0$ is transversely isotopic to a quasipositive braid with respect to $(B',\pi')$, so is $K_1=K_0 \cup U$.

In case (ii-b), the (unique) critical point of odd index must correspond to a positive singular point of the characteristic foliation. The surface then consists of a collection of annuli with linear open book foliation and a pair of pants with positive hyperbolic open book foliation. In this case, $K_1$ is obtained from $K_0$ by braid isotopy and the addition of a positive half-twist.  Since $K_0$ is transversely isotopic to a quasipositive braid $K_0^+$ with respect to $(B',\pi')$, the transverse Markov theorem implies that $K_0$ and $K_0^+$ have a common positive Markov stabilization. Any positive stabilization of a quasipositive braid is quasipositive, so this implies that $K_0$ can be positively stabilized to a quasipositive braid that, for notational convenience, we will also denote by $K_0^+$. By Lemma~\ref{lem:stab}, we can find a positive stabilization $K_1^+$ of $K_1$ that can be obtained from $K_0^+$ by adding a positive half-twist. The addition of a positive half-twist preserves quasipositivity, so $K_1^+$ is also a quasipositive braid with respect to $(B',\pi')$. We conclude that $K_1$ is transversely isotopic to a quasipositive braid.

In each case, it follows that $K_1$ is transversely isotopic to a quasipositive braid $K_1^+$ with respect to the isotoped open book $(B',\pi')$. As above, we may assume that $K_1^+$ avoids the neighborhoods of the critical points of $f: S \to \rr$ where the ambient isotopy from Theorem~\ref{thm:foln} failed to be a contact isotopy. We may then reverse the isotopy, carrying $(B',\pi')$ back to $(B,\pi)$ and carrying $K_1^+$ by a transverse isotopy to a quasipositive braid with respect to the original open book $(B,\pi)$.
\end{proof}

\begin{rem}\label{rem:factor}
We can sharpen the above result:  If $K_0$ is transversely isotopic to the closure of a quasipositive braid with respect to an underlying abstract open book $(F,\varphi)$ with a given factorization of $\varphi$ into Dehn twists, 
 then $K_1$ is also transversely isotopic to the closure of a quasipositive braid with respect to the same factorization of $\varphi$ into Dehn twists. (See Definition~\ref{def:closure} and the discussion preceding it.) This is clear in the case where the characteristic foliation contains no hyperbolic points, since $K_1$ is then obtained from $K_0$ by isotopy and the addition of disjoint unknots meridional to the binding. 
 
Following the inductive argument above, it suffices to consider the case where $S$ is a pair of pants with hyperbolic open book foliation and $\partial S = K_1 \cup - K_0$.  The proof shows that there is \emph{some} pointed open book $(F,P,\widehat \varphi)$ encoding $K_0$ in $(B',\pi')$  such that $\widehat \varphi$ is isotopic to $\varphi$ and  $K_1$ is obtained by adding a positive half-twist $H_\beta$ to $\widehat\varphi$.   By hypothesis, the braid $K_0$ is transversely isotopic to the closure of a quasipositive braid with respect to the given factorization of $\varphi$ into Dehn twists. As above, this quasipositive braid and $K_0$ have a common positive Markov stabilization, which we can write as the closure of a quasipositive braid $H_{\alpha_1} \circ \cdots \circ H_{\alpha_\ell}$ with respect to the given factorization of $\varphi$ into Dehn twists. Mirroring the argument from Lemma~\ref{lem:stab}, we  carry the half-twist $H_\beta$ through the sequence of Markov stabilizations and isotopies to a half-twist along a new arc $\beta'$ so that $K_1$ is transversely isotopic to the closure of the quasipositive braid $H_{\beta'} \circ H_{\alpha_1} \circ \cdots \circ H_{\alpha_\ell}$ with respect to the given factorization of $\varphi$ into Dehn twists.
\end{rem}

\section{Holomorphic curves in Stein surfaces}\label{sec:stein}

 A complex manifold is \emph{Stein} if it admits a proper biholomorphic embedding as a closed subset of some $\cc^n$. We take particular interest in \emph{Stein surfaces}, which are Stein manifolds of complex dimension two, and \emph{Stein cobordisms} between contact manifolds.  Every Stein manifold $X$ admits a so-called \emph{$J$-convex} function $\rho: X\to [0,\infty)$, where $J$ denotes the induced almost-complex structure on $TX$. In our setting, the condition of $J$-convexity is equivalent to requiring that each regular level set $\rho^{-1}(c)$ has a natural contact structure defined by the field of complex tangent lines. This contact structure can also be described as the kernel of the 1-form given by restricting $\eta_\rho=-d^\cc \rho$ to the level set $\rho^{-1}(c)$, where $d^\cc=J^* d$.  
For a thorough exposition of the theory of Stein manifolds, we refer the reader to \cite{ec:book}.

Our proof of Theorem~\ref{thm:qp} is motivated by a simple observation (cf~\cite{bo:qp,fiedler}):   Given a smooth holomorphic curve $\Sigma$ in a Stein surface $X$ and a smooth convex hypersurface $Y \subset X$, consider the intersection $\Sigma \cap Y$. If this intersection is transverse, then $\Sigma \cap Y$ forms a smooth link that is transverse to the natural contact structure $\xi$ on $Y$. Indeed, since $T\Sigma$ and $\xi$ consist of complex lines in $TX$, their intersection at any point has real dimension zero or two. It follows that $\Sigma\cap Y$ is transverse to the contact structure except for those points $p$ at which $\Sigma$ is tangent to $Y$, where $T_p \Sigma= \xi_p \subset T_pY$. 

This suggests a more flexible class of surfaces in Stein domains:

\begin{defn}[{cf~\cite{bo:qp}}]\label{def:ascending}
A smoothly embedded, oriented surface $\Sigma\subset X$ is \emph{ascending} if $\Sigma$ contains no critical points of $\rho|_X$, the restriction $\rho|_\Sigma$ is a Morse function, and, except at its critical points, each level set  $\rho|_\Sigma^{-1}(c)$ is positively transverse to the contact structure on $\rho|_X^{-1}(c)$. 
\end{defn}

Equivalently, to better mirror the original definition from \cite{bo:qp}, one can replace the third condition with the inequality $(-d\rho \wedge d^\cc \rho)|_\Sigma >0 $ at regular points of $\rho|_\Sigma$. In Section~\ref{subsec:ascending}, we present an equivalent characterization of ascending surfaces phrased in terms of characteristic foliations on surfaces in Stein domains, which we define precisely in Definition~\ref{def:stein-char}. We will show that critical points of $\rho|_\Sigma$ have a sign induced by the complex structure on $X$, and complex curves are examples of ascending surfaces whose critical points are all positive (Proposition~\ref{prop:pseudo}). Thus Theorem~\ref{thm:qp} follows from a more general result: 

\begin{thm}\label{thm:full}
If $\Sigma$ is an ascending surface with positive critical points that is properly embedded in a Stein domain $X$, then $K=\partial \Sigma$ is transversely isotopic to a quasipositive braid with respect to a (positive allowable) open book for $\partial X$. 
\end{thm}

We also prove an analogous result for Stein cobordisms in Section~\ref{subsec:full}. In the case of ascending surfaces in the symplectization of $(Y,\xi)$, we have additional control over the choice of compatible open book:

\begin{thm}\label{thm:prod}
Let $\Sigma$ be an ascending surface with positive critical points in a compact piece $Y \times [0,1]$ of the symplectization of $(Y,\xi)$. If $K_0=\Sigma \cap Y_{0}$ is transversely isotopic to the closure of a quasipositive braid with respect to an open book for $(Y,\xi)$, then $K_1=\Sigma \cap Y_{1}$ is transversely isotopic to the closure of a quasipositive braid with respect to the same open book. 
In particular, if $K_0$ is empty, then $K_1$ is transversely isotopic to a quasipositive braid with respect to any open book supporting $(Y,\xi)$.
\end{thm}

Viewing complex curves in $\cc^2$ as lying in the symplectization of $(S^3,\xi_{\st})$, this recovers the theorem of Boileau and Orevkov:

\begin{cor}[Boileau-Orevkov]\label{cor:bo}
The boundary of any smooth, properly embedded complex curve in $B^4 \subset \cc^2$ is quasipositive (with respect to the standard open book).
\end{cor}

The rest of this section is organized as follows: We discuss properties and examples of ascending surfaces in Section~\ref{subsec:ascending}. Then we prove Theorem~\ref{thm:prod} and Corollary~\ref{cor:bo} in Section~\ref{subsec:prod} and Theorem~\ref{thm:full} in Section~\ref{subsec:full}.

\subsection{Smooth and ascending surfaces in Stein manifolds} \label{subsec:ascending} 

We begin by extending the notion of a characteristic foliation to surfaces in a Stein domain $(X,\rho)$ equipped with a 1-form $\eta_\rho = - d^\cc \rho$ as defined above.

\begin{defn}\label{def:stein-char}
Let $\Sigma$ denote a smooth, oriented surface in $(X,\rho)$, and let $\Omega$ be an area form on $\Sigma$. The \emph{characteristic foliation} on $\Sigma$ is the singular foliation directed by the vector field $v$ satisfying $\iota_v \Omega= \eta_\rho |_\Sigma$.
\end{defn}

As with Definition~\ref{def:char}, the characteristic foliation defined above is independent of the choice of area form $\Omega$. We use characteristic foliations to provide an alternative characterization of ascending surfaces:

\begin{lem}\label{lem:stein-char}
A smooth, oriented surface $\Sigma$ in a Stein domain $(X,\rho)$ satisfies $(-d\rho \wedge d^\cc \rho)|_\Sigma >0$ at regular points of $\rho|_\Sigma$ if and only if the characteristic foliation on $\Sigma$ is weakly gradient-like for $\rho|_\Sigma$.
\end{lem}

\begin{proof}For convenience, we write $\eta$ for $\eta_\rho=-d^\cc \rho$.  Fix a positive area form $\Omega$ on $\Sigma$, and define vector fields $v$ and $w$ by the equations $\iota_v \Omega = \eta|_\Sigma$ and $\iota_w \Omega = - d \rho|_\Sigma$.  We will show that, at a regular point $p$ of $\rho|_\Sigma$, we have $d\rho|_\Sigma(v)>0$ if and only if $(d \rho \wedge \eta)|_\Sigma>0$. To this end, we observe
\begin{align}
 \label{eqn:grad}
d \rho|_\Sigma(v) &= -\Omega(w,v)=\Omega(v,w) \\
(d\rho \wedge \eta)|_\Sigma (v,w) &= d\rho|_\Sigma(v) \cdot \eta(w) -0\cdot 0= (\Omega(v,w))^2 \geq0. \label{eqn:diff}
\end{align}
If $d\rho|_\Sigma(v)>0$ at $p$, then $\Omega_p(v,w)>0$ by \eqref{eqn:grad} and thus $(d\rho \wedge \eta)|_\Sigma(v,w)>0$ at $p$ by \eqref{eqn:diff}. It follows that $(d\rho \wedge \eta)|_\Sigma$ is a positive multiple of  $\Omega$ at $p$, ie~$(d\rho \wedge \eta)|_\Sigma>0$. Conversely, $(d\rho \wedge \eta)|_\Sigma > 0$ at $p$ implies that $(d\rho \wedge \eta)|_\Sigma (v,w)$ and $\Omega(v,w)$ have the same sign at $p$. Since the former is positive, so is the latter, which equals $d\rho|_\Sigma(v)$.
\end{proof}

We now provide examples of ascending surfaces via two constructions that will be essential to the proof of Theorem~\ref{thm:symp}. Recall that a contact three-manifold $(Y,\xi)$ has a \emph{symplectization} $\Symp(Y,\xi)$, defined to be  $Y \times \rr$ with the symplectic form $\omega=d(e^t \alpha)$, where $\xi=\ker \alpha$.  Definition~\ref{def:ascending} naturally extends to this setting using the projection $\rho: \Symp(Y,\xi) \to \rr$ and a canonical ``cylindrical'' almost-complex structure inducing the contact structure $\xi$ on each level set $Y_t=\rho^{-1}(t)$.  Alternatively, the contact structure on $Y_t$ is induced by the 1-form $\eta$ obtained as the pullback of $\alpha$ under the projection $\Symp(Y,\xi) \to Y$.

\begin{ex}[Push-in of a convex surface] \label{ex:convex}

Let $S \subset (Y,\xi)$ be a smoothly embedded surface whose boundary $\partial S$ is a positively transverse link and whose characteristic foliation is Morse-Smale and contains no regular closed leaves. Choose a vector field $v$ directing the characteristic foliation. By Theorem B of \cite{smale}, there is a Morse function $f: S\to \rr$ such that $df(v)>0$ at nonzero points of $v$ and $f(S) \subset (-\infty,0]$ with $f|_{\partial S}\equiv 0$. Viewing $Y \times (-\infty,0]$ as part of the symplectization of $(Y,\xi)$, we can define a surface $\Sigma \subset Y \times (-\infty,0]$ to be the graph of $f$ over $S\subset Y$, i.e.~the image of the map $S \to Y \times (-\infty,0]$ defined by $x\mapsto (x,f(x))$. Letting $\rho$ denote the projection $Y \times (-\infty,0] \to (-\infty,0]$, we see that $\rho|_\Sigma$ pulls back to the function $f$ on $S$.  The requirement $df(v)>0$ at regular points implies that regular level sets of $f$ are positively transverse links, so the same is true of the regular level sets of $\rho|_\Sigma$. It follows that $\Sigma$ is an ascending surface in $Y \times (-\infty,0] \subset \Symp(Y,\xi)$. \hfill $\diamond$
\end{ex}

\begin{ex}[Push-in of an immersed Bennequin surface] \label{ex:bennequin} Suppose a transverse braid $K$ in $(Y,\xi)$ is encoded by a pointed open book $(F,P,\widehat \varphi)$ where the monodromy can be factored in the form
$$\widehat \varphi = H_{\alpha_1}^{\epsilon_1} \circ \cdots \circ H_{\alpha_\ell}^{\epsilon_\ell} \circ \varphi$$
where $\varphi$ fixes a collar neighborhood $N$ of $\partial F$ containing the marked points and $\epsilon_j = \pm 1$.  Fix a set of disjoint arcs in $N$ joining the marked points to $\partial F$ and let $S_0 \subset Y$ be the disjoint union of disks swept out by this family of arcs; see Figure~\ref{fig:bennequin}(b).  Now attach 1-handles $h_j$ to $S_0$ as follows: The core of $h_j$ will be a copy of the arc $\alpha_j$ lying in the page $F_{\theta_j}$ for $\theta_j = 2\pi j/(\ell+1)$.  We can extend this core to a once-twisted band $h_j$ whose sign agrees with that of $\epsilon_j$; see Figure~\ref{fig:bennequin}(c). 
 In general, the twisted bands may intersect the interiors of the disks $S_0$ transversely along embedded ``ribbon'' arcs. Let $S$ denote  the resulting  immersed ``Bennequin'' surface for $K$.
 
\begin{figure}\center
\def\svgwidth{400pt} 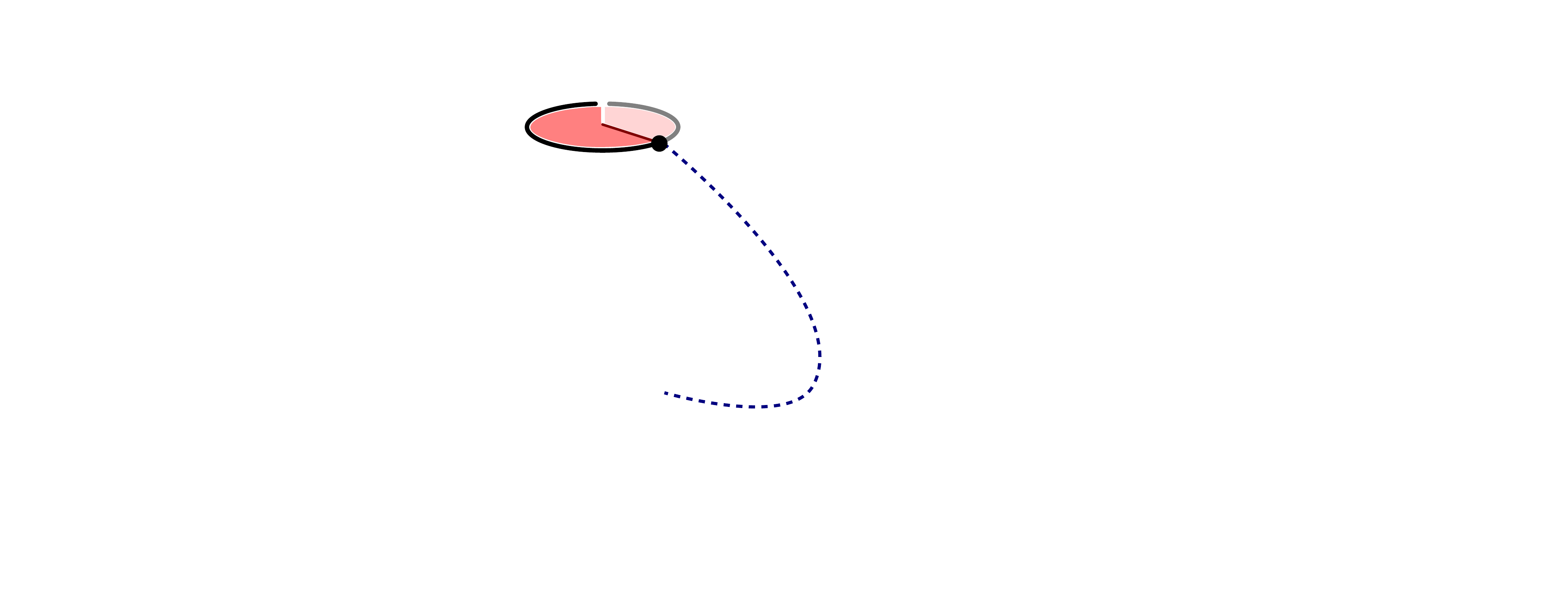 
\caption{Building a Bennequin surface.} \label{fig:bennequin}
\end{figure}

The open book foliation on the immersed surface $S$ is Morse-Smale and has positive elliptic points along the binding and  hyperbolic points along the twisted bands. Using \cite[Theorem 2.21]{ito-kawamuro:open}, we can isotope $\xi$ (away from $\partial S$) to another contact structure $\xi'$ so that the characteristic foliation on $S$ is conjugate to the open book foliation. By Gray's theorem, there is an induced isotopy of $S$ rel boundary so that $\xi$ itself induces the desired characteristic foliation.  Now let $f: S \to \rr$ be a Morse function that has a gradient-like vector field that directs the characteristic foliation. Since any self-intersections of $S$ occur between a disk in $S_0 \subset S$ and part of a twisted band lying in $S \setminus S_0$, we can modify $f$ on the disks (in particular, decrease its value) so that the level sets of $f$ are all embedded. As in Example~\ref{ex:convex}, the graph of $f$ is an embedded ascending surface in $\Symp(Y,\xi)$.  \hfill $\diamond$
\end{ex}

We pause to highlight a corollary. By Example~\ref{ex:bennequin}, every strongly quasipositive transverse link $K$ in $(Y,\xi)$ bounds a Seifert surface whose characteristic foliation is Morse-Smale with positive singularities. By a signed count of elliptic and hyperbolic points in the characteristic or open book foliations on any such surface, we obtain:

\begin{cor}
The Eliashberg-Bennequin inequality is sharp for every strongly quasipositive link in a tight contact three-manifold. \hfill $\square$
\end{cor}

We now consider some general properties of ascending surfaces $\Sigma \subset X$. Recall that a point $p$ in a smooth surface $\Sigma \subset X$ is said to be \emph{complex} if $T_p \Sigma$ is a complex line in $T_p X$. The following proposition extends an observation made by Boileau and Orevkov for ascending surfaces in $\cc^2$.

\begin{prop}\label{prop:complex}
Every critical point of an ascending surface is complex.
\end{prop}

\begin{proof}An ascending surface $\Sigma \subset X$ must be tangent to a level set of $\rho$ at any critical point $p$ of $\rho|_\Sigma$, so it suffices to show that $T_p \Sigma$ lies in the kernel of $\eta=-d^\cc \rho$. Consider the equivalent characterization from the remark following Definition~\ref{def:ascending}: For any area form $\Omega$ on $\Sigma$, we can write
$$(d\rho \wedge \eta)|_\Sigma = f \Omega$$
for some non-negative function $f$ on $\Sigma$ that vanishes only at critical points of $\rho|_\Sigma$.  Since $f$ is non-negative, it vanishes to even order at critical points, hence the same is true for the 2-form $(d\rho \wedge \eta)|_\Sigma$. But since $\rho|_\Sigma$ is Morse, it vanishes precisely to second order at its critical points. This implies that $d\rho|_\Sigma$ vanishes precisely to first order. Since $(d\rho \wedge \eta)|_\Sigma$ vanishes to at least second order, we conclude that $\eta|_\Sigma$ must also vanish at critical points.
\end{proof}

Given an oriented smooth surface $\Sigma \subset X$, every complex point $p\in \Sigma$ can be described as \emph{positive} or \emph{negative} according to whether the intrinsic orientation on  $T_p \Sigma$  agrees or disagrees, respectively, with the complex orientation. In general, we can further classify complex points as \emph{elliptic} or \emph{hyperbolic}; see \cite[Definition~8.3.3]{os:surgery} for the general definition. In the case where $\Sigma$ is an ascending surface, it suffices to define a critical complex point of $\Sigma$ to be \emph{elliptic} if the Morse index of  $\rho|_\Sigma$  is even and \emph{hyperbolic} otherwise. (Note that the \emph{sign} of a critical complex point $p \in \Sigma$ coincides with its sign as a singular point of the characteristic foliation defined in Definition~\ref{def:stein-char}.)

The following proposition shows that Theorem~\ref{thm:qp} is a special case of Theorem~\ref{thm:full}.

\begin{prop}\label{prop:pseudo}
With respect to a generic $J$-convex function, a smooth pseudoholomorphic curve in a Stein surface $(X,J)$ is ascending with positive critical points.
\end{prop}

\begin{proof}
Recall that the subspace of $J$-convex functions $\rho$ on $X$ is open with respect to the $C^2$-topology on the space of all real-valued functions on $X$. Given a pseudoholomorphic curve $\Sigma \subset X$, a generic choice of $\rho$ will ensure that the functions $\rho|_X$ and  $\rho|_\Sigma$ are Morse and have disjoint sets of critical points. As discussed above, all regular level sets of $\rho|_\Sigma$ are transverse links in $\rho^{-1}(t)$. And since $\Sigma$ is pseudoholomorphic, \emph{every} point is complex and positive, including the critical points.
 \end{proof}

\subsection{Ascending surfaces in the symplectization} \label{subsec:prod}

Suppose that $\Sigma$ is a smooth, oriented surface in $\Symp(Y,\xi)$. As discussed in Section~\ref{subsec:ascending}, the 1-form $\eta$ used to define the characteristic foliation on $\Sigma$ is the pullback of a contact form for $\xi$ under the projection $\Symp(Y,\xi) \to Y$. It follows that the characteristic foliation on $\Sigma$ induced by $\eta$ coincides with the characteristic foliation on $\Sigma$ induced by the corresponding map $\Sigma \to Y$.

\begin{proof}[Proof of Theorem~\ref{thm:prod}] Let  $\Sigma$ be an ascending surface with positive critical points in a compact piece $Y \times [0,1]$ of the symplectization of $(Y,\xi)$ with the natural Morse function $\rho$ given by the second-coordinate projection. Suppose that  $K_0=\Sigma \cap Y_0$ is transversely isotopic to a quasipositive braid in a given open book $(B,\pi)$ for $(Y,\xi)$. We wish to show that $K_1=\Sigma \cap Y_{1}$ is also transversely isotopic to a quasipositive braid in $(B,\pi)$. We will   reduce to the three-dimensional setting and apply Theorem~\ref{thm:embedded}.

We may view $K_0$ and $K_1$ as transverse links in $(Y,\xi)$ via the projection $\Symp(Y,\xi) \to Y$. By a slight perturbation of $\Sigma$ (through ascending surfaces), we may assume that $K_0$ and $K_1$ project disjointly to $(Y,\xi)$. These transverse links cobound the singular surface $\phi(\Sigma) \subset Y$ obtained via the composition $\phi:\Sigma \hookrightarrow Y \times [0,1] \to Y$. For technical reasons (see Remark~\ref{rem:embedded-nonsing}), we perturb $\Sigma$ so that the critical level sets of $\rho|_\Sigma$ have disjoint images under $\phi$:  For each critical level set $\rho|_\Sigma^{-1}(t)$, we may perform a small contact isotopy of $Y_t$ that ensures the image of $\rho|_\Sigma^{-1}(t)$ under $\phi$ is disjoint from the image of each other critical level set of $\rho|_\Sigma$. This naturally extends to an isotopy of $Y \times [0,1]$ supported near $Y_t$ that restricts to a contact isotopy on each level hypersurface (and thus preserves the fact that $\Sigma$ is ascending). Repeating this near every critical level set of $\rho|_\Sigma$ induces the desired isotopy of $\Sigma$. After this perturbation, the critical level sets of $\rho|_\Sigma$ are mapped disjointly by $\phi$ into $Y$.

By Lemma~\ref{lem:nonsing}, the map $\phi$ is nonsingular with respect to $\rho|_\Sigma : \Sigma \to \rr$. By Lemma~\ref{lem:stein-char}, the characteristic foliation on $\Sigma$ is weakly gradient-like with respect to $\rho|_\Sigma$ and contains only positive singular points; see also Proposition~\ref{prop:complex} and the ensuing discussion of signs of critical points. We may now apply Theorem~\ref{thm:embedded} (via Remark~\ref{rem:embedded-nonsing}) to conclude that $K_1$ is transversely isotopic to a quasipositive braid with respect to the chosen open book $(B,\pi)$.
\end{proof}

As a corollary, we recover the theorem of Boileau and Orevkov:

\begin{proof}[Proof of Corollary~\ref{cor:bo}]
Let $\Sigma \subset \cc^2$ be a smooth complex curve meeting the unit sphere $S_1^3 \subset \cc^2$ transversely in a link $K$. After a small translation, we may assume that $\Sigma$ lies in $\cc^2 \setminus \{0\} \approx \Symp(S^3,\xi_\st)$ and that the radial distance function restricts to a Morse function on $\Sigma$. (Note that any sufficiently small translation of $\Sigma$ induces a transverse isotopy of the link $K=\Sigma \cap S_1^3$.) By Proposition~\ref{prop:pseudo}, $\Sigma$ is an ascending surface with positive critical points in $\Symp(S^3,\xi_\st)$. Applying Theorem~\ref{thm:prod}, we see that $\Sigma \cap S^3_1$ is quasipositive with respect to any open book supporting $(S^3,\xi_\st)$, including the standard open book $(D^2,\id)$.
\end{proof}

\subsection{The general case}
\label{subsec:full}

Now we consider ascending surfaces in arbitrary Stein domains. By definition, an ascending surface $\Sigma$ in $X$ is sensitive only to the function $\rho$ and the contact structure on level sets of $\rho$.  Through this lens, any noncritical region of a Stein cobordism looks like the symplectization of a regular level set.

\begin{lem}\label{lem:vertical}
Suppose $(X,\rho)$ is a compact Stein cobordism with no critical points and let $(Y,\xi)$ be a regular contact hypersurface $\rho^{-1}(c)$. Then there is a diffeomorphism $$\psi: Y \times [a,b] \subset \Symp(Y,\xi) \to X$$
such that $\psi^* \rho$ is the second-coordinate projection and $\psi$ restricts to a contactomorphism between the level sets of $\rho$ and $\psi^*\rho$.
\end{lem}

To prove Lemma~\ref{lem:vertical}, we take advantage of the \emph{Weinstein structure} on a Stein manifold, which is a triple $(\omega,\rho,v)$ where $\omega$ is a symplectic form, $\rho: W \to \rr$ is an exhausting generalized Morse function, and $v$ is a complete vector field which is Liouville for $\omega$ and gradient-like for $\rho$; see \cite{ec:book} for more detail. On a Stein manifold $X$ with $J$-convex function $\rho: X \to \rr$, we define a Weinstein structure where $\omega_\rho = d\eta_\rho$ (with $\eta_\rho=-d^\cc \rho$ as above) and where $v_\rho$ is gradient vector field for $\rho$ with respect to the Riemannian metric $g_\rho(\cdot,\cdot) =\omega_\rho(\cdot, J \cdot)$.

\begin{proof}[Proof of Lemma~\ref{lem:vertical}] For notational convenience, we assume $c=0$. As a first step, we construct a diffeomorphism $Y \times [a,b] \to X$ that pulls $\rho$ back to the second-coordinate projection. Let $v$ be the gradient of $\rho$ defined above. Since $\rho$ has no critical points and $X$ is compact, the function $f=1/d\rho(v)$ is bounded. By the completeness of $v$ and the boundedness of $f$, the multiple $v'=fv$ is also complete. This vector field satisfies $d\rho(v')\equiv 1$ and is complete, so the flow of $v'$ defines a diffeomorphism from $Y \times [a,b]$ to $X$ that pulls $\rho$ back to the second-coordinate projection.

Now let $\zeta_0$ be the vertically invariant plane field on $Y \times [a,b]$ defined by $\xi$ on each level set, and let $\zeta_1$ be the plane field obtained by pulling back the natural contact structure on each level set of $\rho$ in $X$ under the diffeomorphism constructed above.  We will produce an isotopy of $Y \times [a,b]$ preserving level sets that carries $\zeta_0$ to $\zeta_1$. Choose 1-forms $\alpha_0$ and $\alpha_1$ such that $\zeta_i=\ker(dt \wedge \alpha_i)$.   Now define a family of maps $\pi_s$ from $Y \times [a,b]$ to itself by $\pi_s(y,t)=(y,st)$. 
The map $\pi_s$ restricts to a diffeomorphism from $Y \times \{t\}$ to $Y\times \{st\}$ for all $s \in [0,1]$ and $ t \in [a,b]$. The induced family of 1-forms $\alpha_s=\pi_s^* \alpha_1$ on $Y \times [a,b]$ interpolates between $\alpha_0$ and $\alpha_1$ as defined above because $\pi_s$ is the projection $Y \times [a,b] \to Y \times \{0\}$ for $s=0$ and is the identity for $s=1$. It is also easy to see that  $\alpha_s$ restricts to a contact form on each $Y \times \{t\}$. 

From here, the standard Moser's trick argument produces  a vector field $V_s$ such that 
\begin{enumerate}[label=(\roman*)]
\item $V_s$ is tangent to the level sets $Y \times \{t\}$,
\item $V_s$ lies in $\ker \alpha_s$, and
\item  the flow $\psi_s$ of $V_s$ satisfies
$\psi_s^*\alpha_s=\lambda_s \alpha_0$
for some family of smooth positive functions $\lambda_s$ on $Y \times [a,b]$.
\end{enumerate}
The diffeomorphism $\psi_1$ of $Y \times [a,b]$ preserves level sets and is seen to carry $\zeta_0$ to $\zeta_1$: 
$$(\psi_1)_*\zeta_0=(\psi_1)_*\ker(dt \wedge \alpha_0)=(\psi_1)_*\ker(dt \wedge \lambda_1\alpha_0)=\ker(dt \wedge \alpha_1)=\zeta_1.$$
It follows that the composition of $\psi_1$ with the original map $Y \times [a,b] \to X$ provides the desired diffeomorphism.
\end{proof}

For any Stein domain $X$, we may consider a handle decomposition of $X$ defined by a $J$-convex function $\rho: X \to \rr$. The following proposition says that the the 2-handles have Legendrian attaching circles and that the cores of the 2-handles are Lagrangian (ie~half-dimensional submanifolds on which $\omega$ vanishes).

\begin{prop}[{\cite[Lemma 11.13]{ec:book}}]\label{prop:leg}
Let $(X,J,\rho)$ be a Stein four-manifold with its standard Weinstein structure $(\omega=d\eta,\rho,v)$ as defined above.
\begin{enumerate}[label=\normalfont(\alph*)]
\item The stable manifold $V_p^-$ of any critical point $p \in X$ of $\rho$ satisfies $\eta|_{V_p^-}\equiv 0$. In particular, $V_p^-$ is isotropic for the symplectic structure $\omega$ and the intersection $V_p^- \cap \rho^{-1}(c)$ with any regular level set is isotropic for the contact structure induced by $\eta$ on $\rho^{-1}(c)$.

\item Suppose $\rho$ has no critical values in $[a,b]$. Then the image of any isotropic submanifold  $\Lambda^a \subset \rho^{-1}(a)$ under the flow of $v$ intersects $\rho^{-1}(b)$ in another isotropic submanifold $\Lambda^b$.
\end{enumerate}
\end{prop}

We obtain Theorem~\ref{thm:full} as an immediate corollary of the following stronger statement for ascending surfaces in Stein cobordisms: 

\begin{thm}\label{thm:cobordism}
Let $\Sigma$ be an ascending surface with positive critical points in a Stein cobordism $X$ between $(Y_0,\xi_0)$ and $(Y_1,\xi_1)$. If $K_0 = \Sigma \cap Y_0$ is transversely isotopic to a quasipositive braid in a (positive allowable) open book compatible with $\xi_0$, then $K_1=\Sigma \cap Y_1$ is transversely isotopic to a quasipositive braid in a (positive allowable) open book compatible with $\xi_1$.
\end{thm}

\begin{proof}Let $(B_0,\pi_0)$ be a compatible open book for $(Y_0,\xi_0)$ with respect to which $K_0$ is quasipositive. We will show that $K_1$ is quasipositive with respect to a compatible open book $(B_1,\pi_1)$ for $(Y_1,\xi_1)$ obtained from $(B_0,\pi_0)$ by positive Hopf stabilization and the addition of positive Dehn twists (along homologically nontrivial curves) into the monodromy. By subdividing $X$, it suffices to consider the case where $X$ contains at most one critical point. Note that, by the $C^2$-openness of $J$-convex functions, we may choose $\rho$ so that the critical values of $\rho|_X$ and $\rho|_\Sigma$ are all distinct.

If $X$ contains no critical points of $\rho$, then for any value $c$ there's a diffeomorphism from a compact piece of $\Symp(\rho^{-1}(c))$ to $X$ of the form described in Lemma~\ref{lem:vertical}. By that lemma and the definition of an ascending surface, it is easy to see that $\Sigma \subset X$ pulls back to an ascending surface with positive critical points in $\Symp(\rho^{-1}(c))$. The claim then follows from Theorem~\ref{thm:prod}.

Now suppose that $\rho|_X$ has a single critical point, the value of which we will denote by $c$. For any $\epsilon>0$, we can further subdivide $X$ into three pieces: two exterior pieces where $|\rho -c |>\epsilon$ and an interior piece where $|\rho -c| \leq \epsilon$. It suffices to prove the claim for the interior piece. First consider the case where the critical point of $\rho|_X$ has Morse index zero. This is equivalent to considering ascending surfaces that miss the origin in the standard Stein $B^4 \subset \cc^2$, so the ascending surface can be viewed in $\Symp(S^3,\xi_\st)$,  where the desired claim again follows from Theorem~\ref{thm:prod}.

Next consider the case where the critical point of $\rho|_X$ has Morse index one. By choosing $\epsilon$ sufficiently small, we can assume that $\rho|_\Sigma$ has no critical points in $X$. The level set $Y_1=\rho^{-1}(c+\epsilon)$ is obtained from $Y_0=\rho^{-1}(c-\epsilon)$ by removing a pair of small Darboux balls and attaching a copy of $S^2 \times I$ with the appropriate contact structure. After positive Hopf stabilization, we can assume that $Y_0$ is equipped with an open book $(B_0,\pi_0)$ such that (1) the attaching region $A$ for the 1-handle lies along the binding $B_0$ and (2) the transverse link $K_0=\Sigma \cap Y_0$ is transversely isotopic to a quasipositive braid $K_0^+$ with respect to $(B_0,\pi_0)$. (Note that positive Hopf stabilization of a pointed open book  preserves quasipositivity.) Moreover, since the attaching region $A$ is a pair of small Darboux balls,  we may assume that the transverse isotopy between $K_0$ and $K_0^+$ is performed in the complement of $A$. Now we construct a natural open book on $Y_1$ that extends the open book on $Y_0$: By an argument similar to the one used in the proof of Lemma~\ref{lem:vertical}, there is an embedding
$$\Psi:(Y_0 \setminus A) \times [c-\epsilon,c+\epsilon] \hookrightarrow X$$
pulling $\rho$ back to the second-coordinate projection and restricting to a contact embedding on each level set. Since $\Sigma$ avoids the critical points of $\rho|_X$, we can also assume that $\Sigma$ lies in the image of this embedding. Define an embedding $\psi: Y_0 \setminus A \hookrightarrow Y_1$ by sending $y \in Y_0\setminus A$ to $\Psi(y,c+\epsilon)$ in $Y_1$. We can carry the partial open book from $Y_0 \setminus A$ to $Y_1$ using $\psi$ and then extend it in the natural way over the 1-handle in $Y_1$ to define a compatible open book $(B_1,\pi_1)$. Abstractly, this can be described as follows:  if $(B_0,\pi_0)$ is encoded by $(F_0,\varphi_0)$, then $(B_1,\pi_1)$ is encoded by $(F_1,\varphi_1)$ where $F_1$ is the union of $F_0$ and an oriented 1-handle and $\varphi_1$ is the extension of $\varphi_0$ to $F_1$ via the identity on the 1-handle. 

Now the quasipositive braid $K_0^+$ with respect to $(B_0,\pi_0)$ in $Y_0 \setminus A$ induces a quasipositive braid $K_1^+=\psi(K_0^+)$ with respect to $(B_1,\pi_1)$ in $\psi(Y_0 \setminus A) \subset Y_1$. Since the transverse isotopy between $K_0$ and $K_0^+$ lies in the domain of $\psi$, the transverse link $\psi(K_0)$ in $Y_1$ is transversely isotopic to $K_1^+$ in $Y_1$. And since $\Sigma$ lies in the image of $\Psi$ and we have assumed $\rho|_\Sigma$ has no critical points in $X$, we can use $\Psi$ and $\Sigma$ to realize a transverse isotopy between $\psi(K_0)$ and $K_1$ in $Y_1$. It follows that $K_1$ is transversely isotopic to the quasipositive braid $K_1^+$ with respect to $(B_1,\pi_1)$.

Finally, we must consider the case where the critical point of $\rho|_X$ has Morse index two. In this case, the argument is nearly identical to the preceding case. The primary difference is that the attaching region $A$ of the Stein 2-handle is a neighborhood of a Legendrian knot $\Lambda$ in $(Y_0,\xi_0)$; see Proposition~\ref{prop:leg}. However, any open book compatible with $\xi_0$ has a positive Hopf stabilization containing the Legendrian knot $\Lambda$ in a page; see \cite[Corollary 4.23]{etnyre:obd}. Given an open book with the desired form, we now note that the contact three-manifold $(Y_1,\xi_1)$ is obtained from $(Y_0,\xi_0)$ by Legendrian surgery along $\Lambda$. As above, we can construct a family of contact embeddings $Y_0 \setminus A \hookrightarrow \rho^{-1}(t)$ for $t$ near $c$ and produce the desired open book $(B_1,\pi_1)$ from $(B_0,\pi_0)$. Abstractly, this corresponds to inserting a positive Dehn twist along the representative of $\Lambda$ in the page of the open book. Note that, since any contact form $\alpha$ for $\xi_0$ induces an area form $d\alpha$ on each page of $(B_0,\pi_0)$,  the Legendrian curve $\Lambda$ must be homologically nontrivial by Stokes' theorem.

It remains to show that $K_1$ is transversely isotopic to a quasipositive braid with respect to $(B_1,\pi_1)$. 
The family of contact embeddings $Y_0 \setminus A \hookrightarrow  \rho^{-1}(t)$ also lets us pull $\Sigma$ back to an ascending surface in $\Symp(Y_0 \setminus A, \xi_0)$ with no critical points. By the transverse Alexander and Markov theorems, we can find braid representatives of $K_0$ and $K_1$ encoded by pointed open books that are related by positive Markov stabilization and destabilization and the addition of a positive Dehn twist.  By an argument analogous to the proof of Lemma~\ref{lem:stab}, it is easy to show that these operations preserve the property of being transversely isotopic to a quasipositive braid.

In the case where $K_0= \Sigma \cap Y_0$ is transversely isotopic to the closure of a quasipositive braid with respect to a positive allowable open book, it is clear that the same is true of $K_1$ because we have modified the original underlying open book only by adding 1-handles, inserting positive Dehn twists along homologically nontrivial curves in the page, and performing positive Hopf stabilizations (which also introduce Dehn twists along homologically nontrivial curves). See also Remark~\ref{rem:factor}.
\end{proof}

Observe that Theorem~\ref{thm:full} follows by taking $K_0$ to be empty and viewing a Stein domain $X$ as a cobordism from the empty set to $\partial X$. Theorem~\ref{thm:qp}  follows in turn from Theorem~\ref{thm:full} and Proposition~\ref{prop:pseudo}.

\section{Constructions of symplectic surfaces}\label{sec:con}

In this final section, we prove that every quasipositive braid in a positive allowable open book bounds a symplectic surface in some Stein filling of the associated contact three-manifold. Our strategy will be to construct ascending surfaces with positive critical points and then perform an isotopy to make the surfaces symplectic. This requires a lemma adapted from \cite{bo:qp}.

\begin{lem}\label{lem:stretch}
Let $\Sigma$ be an ascending cobordism with positive critical points in $\Symp(Y,\xi)$ between transverse links $K_0$ and $K_1$ in $(Y,\xi)$. Then $\Sigma$ can be made symplectic through a compactly supported isotopy of ascending surfaces.

\end{lem}

\begin{proof}
For any constant $a \geq 1$, define a diffeomorphism  $\psi_a$ of $Y \times \rr$ by $\psi_a(y,t)=(y,at).$
 Assume that $\Sigma$ is cylindrical outside of $[T_-,T_+]$ and let $\Sigma'=\Sigma \cap [T_-,T_+]$. We will first show that for $a \gg 1$, we have $\omega|_{\psi_a(\Sigma')}>0$, i.e. $\psi_a^*(\omega)|_{\Sigma'}>0$. An elementary calculation gives
\begin{align*}
\psi_a^*(\omega)
=a e^{at} dt \wedge \alpha +e^{at} d\alpha
= e^{(a-1)t}\left( (a-1)e^t\, dt \wedge \alpha + \omega\right).
\end{align*}
Recall that $(dt \wedge \alpha)|_\Sigma \geq 0$, with equality if and only if the point in $\Sigma$ is a singular point of $\rho|_\Sigma$. In terms of an area form $\Omega$ for $\Sigma$, this means we can write $e^t(dt \wedge \alpha)|_\Sigma$ as $f \cdot \Omega$ with $f \geq 0$. Similarly, we can write $\omega=g \cdot \Omega$. Note that $g>0$ on a neighborhood of each critical point of $\rho|_\Sigma$ because $\Sigma$ is ascending with positive critical points. Using the compactness of $\Sigma'$, it follows that $\psi_a^*(\omega)|_{\Sigma'}>0$ for $a \gg 1$ and thus $\psi_a(\Sigma')$ is symplectic. Outside of $\psi_a(\Sigma')$, the surface $\psi_a(\Sigma)$ is cylindrical over transverse links. It follows that $\psi_a(\Sigma)$ is a symplectic cobordism between $K_0$ and $K_1$. Moreover, it is easy to check that $\psi_a(\Sigma)$ remains ascending with positive critical points.
\end{proof}

\begin{rem} 
If $(Y,\xi)$ is symplectically fillable and $K_0=\emptyset$, the relative symplectic Thom conjecture implies that the surfaces constructed are genus-minimizing.
\end{rem}

The desired surfaces can now be constructed using  Example~\ref{ex:bennequin} and Lemma~\ref{lem:stretch}.

\begin{proof}[Proof of Theorem~\ref{thm:symp}]
 Let $(Y,\xi)$ denote the contact three-manifold in question. By hypothesis, we have a transverse braid $K$ in $Y$ encoded by a Stein quasipositive pointed open book $(F,P, \widehat \psi \circ  \varphi)$ where $\widehat \psi$ is a product of positive half-twists and $\varphi=D_{\gamma_1} \circ \cdots \circ D_{\gamma_m}$ is a product of positive Dehn twists along homologically nontrivial simple closed curves $\gamma_j$ in $F$. The modified pointed open book $(F,P,\widehat \psi \circ \id)$ is a quasipositive braid $K_0$ in $Y_0=\#^k(S^1 \times S^2)$ for $k=2g(F)+|\partial F|-1$. The associated contact structure $\xi_0$ is the unique Stein-fillable contact structure obtained from viewing $Y_0$ as the boundary of the subcritical Stein domain $X_0 = \natural^k(S^1 \times B^3)$. As constructed in Example~\ref{ex:bennequin}, we can find an ascending surface $\Sigma$ with positive critical points in a compact piece $Y_0 \times (-\epsilon,0]$ of the symplectization $\Symp(Y_0,\xi_0)$. By applying Lemma~\ref{lem:stretch} and allowing $\epsilon$ to increase, we can isotope this surface to be symplectic. After stretching a collar neighborhood of $\partial X_0$, the surface $\Sigma$ in $Y_0 \times (-\epsilon,0]$ embeds properly and symplectically into $X_0$. By construction, the boundary of $\Sigma$ is the transverse link $K_0$ encoded by $(F,P,\widehat \psi \circ \id)$.

Using the Legendrian Realization Principle (see \cite{yamamoto:overtwisted} for this special case), we can realize the homologically nontrivial curves $\gamma_j$ as Legendrian knots $\Lambda_j$ lying on  distinct pages $F_{\theta_j}$ in $(Y_0,\xi_0)$ for increasing $\theta_j\in (\pi,2\pi)$; here we scale $\theta$ so that all nontrivial braiding occurs while $\theta$ lies in $(0,\pi)$. Performing surgery along each $\Lambda_j$ with framing $-1$ (relative to the page framing) corresponds to modifying the pointed open book by introducing a positive Dehn twist along each $\gamma_j$ in order of increasing $\theta_j$. The result reproduces the original pointed open book $(F,P,\widehat \psi \circ \varphi)$. These surgeries can be realized by attaching Stein 2-handles along the Legendrian knots $\Lambda_j$ to produce a Stein filling $X$ of $(Y,\xi)$. Moreover, since the handles are attached away from the boundary of $\Sigma$, the symplectic structure near $\Sigma$ is unchanged. It follows that $\Sigma$ is a symplectic surface in $X$ whose boundary is the original link $K$.
\end{proof}

Combining this construction with Theorem~\ref{thm:full} yields the following characterization:

\begin{cor}\label{cor:char}
A transverse link in $(Y,\xi)$ is quasipositive with respect to a positive allowable open book if and only if it bounds an ascending surface with positive critical points in some Stein filling of $(Y,\xi)$. \hfill $\square$
\end{cor}

We end this section with an example that demonstrates the failure of Theorem~\ref{thm:symp} for links in the boundary of a Stein domain that are quasipositive braids with respect to an open book whose monodromy is \emph{not} positive allowable.

\begin{ex}Consider an annulus $A$ with a single marked point $p$ as depicted in Figure~\ref{fig:push}. Letting $P_\delta$ denote the push map  about the core circle through $p$, we note that the pointed open book $(A,\{p\},P_\delta)$ encodes a braid $K$ in $S^1 \times S^2$ with respect to its trivial planar open book. This braid $K$ represents a nontrivial element in $H_1(S^1 \times S^2)$.  Since any Stein filling of $S^1{\times}S^2$ is diffeomorphic to $S^1{\times}B^3$ and $H_2(S^1{\times}B^3,\partial(S^1{\times}B^3))=0$, we see that $K$ cannot bound a surface in any Stein filling. However, as depicted on the right side of Figure~\ref{fig:push}, $K$ can also be represented by a pointed monodromy consisting of a pair of (oppositely-signed) Dehn twists. This braid is quasipositive with respect to an open book monodromy that is \emph{not} positive allowable, though the underlying open book supports the Stein-fillable contact structure on $S^1 \times S^2$. \hfill $\diamond$
\end{ex}

%
%
%

\bibliographystyle{gtart}
\bibliography{biblio}

\begin{thebibliography}{}
\providecommand\bibmarginpar{\leavevmode\marginpar}
\def\urlstyle#1{{\tt #1}}

\bibitem{bvvv:equivalence}
\textbf{J\,A Baldwin}, \textbf{D\,S Vela-Vick}, \textbf{V V\'ertesi},
  \href{http://dx.doi.org/10.2140/gt.2013.17.925} {\emph{On the equivalence of
  {L}egendrian and transverse invariants in knot {F}loer homology}}, Geom.
  Topol. 17 (2013) 925--974

\bibitem{square}
\textbf{I Baykur}, \textbf{J Etnyre}, \textbf{M Hedden}, \textbf{K Kawamuro},
  \textbf{J\,V Horn-Morris}, \emph{Contact and symplectic geometry and the
  mapping class groups}, Official report of the 2nd SQuaRE meeting, American
  Institute of Mathematics (July 2015)

\bibitem{bennequin}
\textbf{D Bennequin}, \emph{Entrelacements et Equations de {P}faff}, Asterisque
  107--108 (1983) 87--161

\bibitem{birman-finkelstein}
\textbf{J\,S Birman}, \textbf{E Finkelstein},
  \href{http://dx.doi.org/10.1142/S0218216598000176} {\emph{Studying surfaces
  via closed braids}}, J. Knot Theory Ramifications 7 (1998) 267--334

\bibitem{b-m:simplicity}
\textbf{J\,S Birman}, \textbf{W\,W Menasco},
  \href{http://dx.doi.org/10.2140/gt.2006.10.1425} {\emph{Stabilization in the
  braid groups. {II}. {T}ransversal simplicity of knots}}, Geom. Topol. 10
  (2006) 1425--1452 (electronic)

\bibitem{bo:qp}
\textbf{M Boileau}, \textbf{S Orevkov}, \emph{Quasi-positivit\'e d'une courbe
  analytique dans une boule pseudo-convexe}, C. R. Acad. Sci. Paris S\'er. I
  Math. 332 (2001) 825--830

\bibitem{ec:book}
\textbf{K Cieliebak}, \textbf{Y Eliashberg},
  \href{http://dx.doi.org/10.1090/coll/059} {\emph{From {S}tein to {W}einstein
  and back}}, volume~59 of \emph{American Mathematical Society Colloquium
  Publications}, American Mathematical Society, Providence, RI (2012)Symplectic
  geometry of affine complex manifolds

\bibitem{etnyre:obd}
\textbf{J\,B Etnyre}, \emph{Lectures on open book decompositions and contact
  structures}, from ``Floer homology, gauge theory, and low-dimensional
  topology'', Clay Math. Proc. 5, Amer. Math. Soc., Providence, RI (2006)
  103--141

\bibitem{hedden-etnyre:email}
\textbf{J\,B Etnyre}, \textbf{M Hedden}, Personal communication (February 28,
  2017)

\bibitem{E-VHM:monoids}
\textbf{J\,B Etnyre}, \textbf{J Van Horn-Morris}, \emph{Monoids in the mapping
  class group}, Geometry and Topology Monographs 19 (December 2015) 319--365

\bibitem{farb-margalit}
\textbf{B Farb}, \textbf{D Margalit}, \emph{A primer on mapping class groups},
  volume~49 of \emph{Princeton Mathematical Series}, Princeton University
  Press, Princeton, NJ (2012)

\bibitem{fiedler}
\textbf{T Fiedler}, \href{https://doi.org/10.1007/BF01393886} {\emph{Complex
  plane curves in the ball}}, Invent. Math. 95 (1989) 479--506

\bibitem{gabai:murasugi}
\textbf{D Gabai}, \href{https://doi.org/10.1090/conm/020/718138} {\emph{The
  {M}urasugi sum is a natural geometric operation}}, from ``Low-dimensional
  topology ({S}an {F}rancisco, {C}alif., 1981)'', Contemp. Math. 20, Amer.
  Math. Soc., Providence, RI (1983)  131--143

\bibitem{gk:caps}
\textbf{S Gadgil}, \textbf{D Kulkarni},
  \href{http://dx.doi.org/10.1007/s12044-016-0278-3} {\emph{Relative symplectic
  caps, 4-genus and fibered knots}}, Proc. Indian Acad. Sci. Math. Sci. 126
  (2016) 261--275

\bibitem{geiges:intro}
\textbf{H Geiges}, \emph{An introduction to contact topology}, volume 109 of
  \emph{Cambridge Studies in Advanced Mathematics}, Cambridge University Press,
  Cambridge (2008)

\bibitem{giroux:geometry}
\textbf{E Giroux}, \emph{G\'eom\'etrie de contact: de la dimension trois vers
  les dimensions sup\'erieures}, from ``Proceedings of the {I}nternational
  {C}ongress of {M}athematicians, {V}ol. {II} ({B}eijing, 2002)'', Higher Ed.
  Press, Beijing (2002)  405--414

\bibitem{hedden:subcritical}
\textbf{M Hedden}, \emph{Knot theory and algebraic curves in subcritical
  {S}tein domains}, In preparation.

\bibitem{hedden:email}
\textbf{M Hedden}, Private communication (November 11, 2016)

\bibitem{hkm:sutured}
\textbf{K Honda}, \textbf{W\,H Kazez}, \textbf{G Mati\'c},
  \href{http://dx.doi.org/10.1007/s00222-008-0173-3} {\emph{The contact
  invariant in sutured {F}loer homology}}, Invent. Math. 176 (2009) 637--676

\bibitem{ito-kawamuro:open}
\textbf{T Ito}, \textbf{K Kawamuro},
  \href{http://dx.doi.org/10.2140/gt.2014.18.1581} {\emph{Open book
  foliation}}, Geom. Topol. 18 (2014) 1581--1634

\bibitem{ik:bennequin}
\textbf{T Ito}, \textbf{K Kawamuro}, \emph{The defect of Bennequin-Eliashberg
  inequality and Bennequin surfaces} (March 2017)Available as arXiv:1703.09322

\bibitem{mcduff:positivity}
\textbf{D McDuff}, \href{http://dx.doi.org/10.1007/978-3-0348-8508-9_7}
  {\emph{Singularities and positivity of intersections of {$J$}-holomorphic
  curves}}, from ``Holomorphic curves in symplectic geometry'', Progr. Math.
  117, Birkh\"{a}user, Basel (1994)  191--215With an appendix by Gang Liu

\bibitem{micallef-white}
\textbf{M\,J Micallef}, \textbf{B White},
  \href{http://dx.doi.org/10.2307/2118627} {\emph{The structure of branch
  points in minimal surfaces and in pseudoholomorphic curves}}, Ann. of Math.
  (2) 141 (1995) 35--85

\bibitem{os:markov}
\textbf{S\,Y Orevkov}, \textbf{V\,V Shevchishin},
  \href{http://dx.doi.org/10.1142/S0218216503002846} {\emph{Markov theorem for
  transversal links}}, J. Knot Theory Ramifications 12 (2003) 905--913

\bibitem{os:surgery}
\textbf{B Ozbagci}, \textbf{A\,I Stipsicz},
  \href{http://dx.doi.org/10.1007/978-3-662-10167-4} {\emph{Surgery on contact
  3-manifolds and {S}tein surfaces}}, volume~13 of \emph{Bolyai Society
  Mathematical Studies}, Springer-Verlag, Berlin; J\'anos Bolyai Mathematical
  Society, Budapest (2004)

\bibitem{pavelescu:braids}
\textbf{E Pavelescu}, \href{http://dx.doi.org/10.2140/pjm.2011.253.475}
  {\emph{Braiding knots in contact 3-manifolds}}, Pacific J. Math. 253 (2011)
  475--487

\bibitem{rudolph:qp-alg}
\textbf{L Rudolph}, \emph{Algebraic functions and closed braids}, Topology 22
  (1983) 191--202

\bibitem{rudolph:kauffman-bound}
\textbf{L Rudolph}, \emph{A congruence between link polynomials}, Math. Proc.
  Cambridge Philos. Soc. 107 (1990) 319--327

\bibitem{rudolph:qp-obstruction}
\textbf{L Rudolph}, \emph{Quasipositivity as an obstruction to sliceness},
  Bull. Amer. Math. Soc. (N.S.) 29 (1993) 51--59

\bibitem{smale}
\textbf{S Smale}, \emph{On gradient dynamical systems}, Ann. of Math. (2) 74
  (1961) 199--206

\bibitem{tw:existence}
\textbf{W\,P Thurston}, \textbf{H\,E Winkelnkemper},
  \href{http://dx.doi.org/10.2307/2040160} {\emph{On the existence of contact
  forms}}, Proc. Amer. Math. Soc. 52 (1975) 345--347

\bibitem{wrinkle}
\textbf{N\,C Wrinkle},
  \href{http://gateway.proquest.com.ezproxy.cul.columbia.edu/openurl?url_ver=Z39.88-2004&rft_val_fmt=info:ofi/fmt:kev:mtx:dissertation&res_dat=xri:pqdiss&rft_dat=xri:pqdiss:3048264}
  {\emph{The {M}arkov theorem for transverse knots}}, ProQuest LLC, Ann Arbor,
  MI (2002)Thesis (Ph.D.)--Columbia University

\bibitem{yamamoto:overtwisted}
\textbf{R Yamamoto}, \href{http://projecteuclid.org/euclid.jmsj/1191591856}
  {\emph{Open books supporting overtwisted contact structures and the
  {S}tallings twist}}, J. Math. Soc. Japan 59 (2007) 751--761

\end{thebibliography}


\begin{thebibliography}{}
\providecommand\bibmarginpar{\leavevmode\marginpar}
\def\urlstyle#1{{\tt #1}}

\bibitem{bcv:ribbons}
\textbf{S Baader}, \textbf{K Cieliebak}, \textbf{T Vogel},
  \href{http://dx.doi.org/10.1142/S0218216509006999} {\emph{Legendrian ribbons
  in overtwisted contact structures}}, J. Knot Theory Ramifications 18 (2009)
  523--529

\bibitem{bi:leg-qp}
\textbf{S Baader}, \textbf{M Ishikawa},
  \href{http://afst.cedram.org/item?id=AFST_2009_6_18_2_285_0}
  {\emph{Legendrian graphs and quasipositive diagrams}}, Ann. Fac. Sci.
  Toulouse Math. (6) 18 (2009) 285--305

\bibitem{bvvv:equivalence}
\textbf{J\,A Baldwin}, \textbf{D\,S Vela-Vick}, \textbf{V V\'ertesi},
  \href{http://dx.doi.org/10.2140/gt.2013.17.925} {\emph{On the equivalence of
  {L}egendrian and transverse invariants in knot {F}loer homology}}, Geom.
  Topol. 17 (2013) 925--974

\bibitem{bennequin}
\textbf{D Bennequin}, \emph{Entrelacements et Equations de {P}faff}, Asterisque
  107--108 (1983) 87--161

\bibitem{birman-finkelstein}
\textbf{J\,S Birman}, \textbf{E Finkelstein},
  \href{http://dx.doi.org/10.1142/S0218216598000176} {\emph{Studying surfaces
  via closed braids}}, J. Knot Theory Ramifications 7 (1998) 267--334

\bibitem{bo:qp}
\textbf{M Boileau}, \textbf{S Orevkov}, \emph{Quasi-positivit\'e d'une courbe
  analytique dans une boule pseudo-convexe}, C. R. Acad. Sci. Paris S\'er. I
  Math. 332 (2001) 825--830

\bibitem{ec:book}
\textbf{K Cieliebak}, \textbf{Y Eliashberg},
  \href{http://dx.doi.org/10.1090/coll/059} {\emph{From {S}tein to {W}einstein
  and back}}, volume~59 of \emph{American Mathematical Society Colloquium
  Publications}, American Mathematical Society, Providence, RI (2012)Symplectic
  geometry of affine complex manifolds

\bibitem{ef:trivial}
\textbf{Y Eliashberg}, \textbf{M Fraser},
  \href{http://proxy.bc.edu/login?url=http://search.ebscohost.com.proxy.bc.edu/login.aspx?direct=true&db=msn&AN=MR2496415&site=ehost-live}
  {\emph{Topologically trivial Legendrian knots.}}, The Journal of Symplectic
  Geometry 7 (2009) 77

\bibitem{etnyre:obd}
\textbf{J\,B Etnyre}, \emph{Lectures on open book decompositions and contact
  structures}, from ``Floer homology, gauge theory, and low-dimensional
  topology'', Clay Math. Proc. 5, Amer. Math. Soc., Providence, RI (2006)
  103--141

\bibitem{hedden-etnyre:email}
\textbf{J\,B Etnyre}, \textbf{M Hedden}, Personal communication (February 28,
  2017)

\bibitem{E-VHM:monoids}
\textbf{J\,B Etnyre}, \textbf{J Van Horn-Morris}, \emph{Monoids in the mapping
  class group}, Geometry and Topology Monographs 19 (December 2015) 319--365

\bibitem{gk:caps}
\textbf{S Gadgil}, \textbf{D Kulkarni},
  \href{http://dx.doi.org/10.1007/s12044-016-0278-3} {\emph{Relative symplectic
  caps, 4-genus and fibered knots}}, Proc. Indian Acad. Sci. Math. Sci. 126
  (2016) 261--275

\bibitem{geiges:intro}
\textbf{H Geiges}, \emph{An introduction to contact topology}, volume 109 of
  \emph{Cambridge Studies in Advanced Mathematics}, Cambridge University Press,
  Cambridge (2008)

\bibitem{giroux:geometry}
\textbf{E Giroux}, \emph{G\'eom\'etrie de contact: de la dimension trois vers
  les dimensions sup\'erieures}, from ``Proceedings of the {I}nternational
  {C}ongress of {M}athematicians, {V}ol. {II} ({B}eijing, 2002)'', Higher Ed.
  Press, Beijing (2002)  405--414

\bibitem{hayden:thesis}
\textbf{K Hayden},
  \href{http://gateway.proquest.com.ezproxy.cul.columbia.edu/openurl?url_ver=Z39.88-2004&rft_val_fmt=info:ofi/fmt:kev:mtx:dissertation&res_dat=xri:pqm&rft_dat=xri:pqdiss:10787187}
  {\emph{Knot {T}heory of {H}olomorphic {C}urves in {S}tein {S}urfaces}},
  ProQuest LLC, Ann Arbor, MI (2018)Thesis (Ph.D.)--Boston College

\bibitem{hedden:subcritical}
\textbf{M Hedden}, \emph{Knot theory and algebraic curves in subcritical
  {S}tein domains}In preparation

\bibitem{hedden:email}
\textbf{M Hedden}, Private communication (November 11, 2016)

\bibitem{ito-kawamuro:open}
\textbf{T Ito}, \textbf{K Kawamuro},
  \href{http://dx.doi.org/10.2140/gt.2014.18.1581} {\emph{Open book
  foliation}}, Geom. Topol. 18 (2014) 1581--1634

\bibitem{ik:bennequin}
\textbf{T Ito}, \textbf{K Kawamuro}, \emph{The defect of Bennequin-Eliashberg
  inequality and Bennequin surfaces} (March 2017)Available as arXiv:1703.09322

\bibitem{os:markov}
\textbf{S\,Y Orevkov}, \textbf{V\,V Shevchishin},
  \href{http://dx.doi.org/10.1142/S0218216503002846} {\emph{Markov theorem for
  transversal links}}, J. Knot Theory Ramifications 12 (2003) 905--913

\bibitem{pavelescu:braids}
\textbf{E Pavelescu}, \href{http://dx.doi.org/10.2140/pjm.2011.253.475}
  {\emph{Braiding knots in contact 3-manifolds}}, Pacific J. Math. 253 (2011)
  475--487

\bibitem{rudolph:qp-alg}
\textbf{L Rudolph}, \emph{Algebraic functions and closed braids}, Topology 22
  (1983) 191--202

\bibitem{rudolph:kauffman-bound}
\textbf{L Rudolph}, \emph{A congruence between link polynomials}, Math. Proc.
  Cambridge Philos. Soc. 107 (1990) 319--327

\bibitem{rudolph:qp-obstruction}
\textbf{L Rudolph}, \emph{Quasipositivity as an obstruction to sliceness},
  Bull. Amer. Math. Soc. (N.S.) 29 (1993) 51--59

\bibitem{smale}
\textbf{S Smale}, \emph{On gradient dynamical systems}, Ann. of Math. (2) 74
  (1961) 199--206

\bibitem{tw:existence}
\textbf{W\,P Thurston}, \textbf{H\,E Winkelnkemper},
  \href{http://dx.doi.org/10.2307/2040160} {\emph{On the existence of contact
  forms}}, Proc. Amer. Math. Soc. 52 (1975) 345--347

\bibitem{wrinkle}
\textbf{N\,C Wrinkle}, \emph{The Markov theorem for transverse knots}
  (2002)Preprint, arXiv:math/0202055v1[math.GT]

\bibitem{yamamoto:overtwisted}
\textbf{R Yamamoto}, \href{http://projecteuclid.org/euclid.jmsj/1191591856}
  {\emph{Open books supporting overtwisted contact structures and the
  {S}tallings twist}}, J. Math. Soc. Japan 59 (2007) 751--761

\end{thebibliography}

%
%
%

\end{document}